\documentclass{amsart}
\usepackage{amssymb}
\usepackage{amsfonts}
\usepackage{amsmath}
\usepackage{url}
\usepackage{color}
\usepackage{setspace}
\usepackage{enumerate}

\usepackage{epsfig}
\usepackage{graphicx}
\usepackage{psfrag}
\newtheorem{theorem}{Theorem}

\newtheorem{proposition}[theorem]{Proposition}
\newtheorem{corollary}[theorem]{Corollary}

\newtheorem{claim}[theorem]{Claim}

\theoremstyle{remark}
\newtheorem{remark}[theorem]{Remark}
\newtheorem{observation}[theorem]{Observation}

\newtheorem{definition}[theorem]{Definition}
\newtheorem{fact}[theorem]{Fact}

\newcommand{\scalarproduct}{\diamond}

\begin{document}

\title{Canonical triangulations of Dehn fillings}
\author[F. Gu\'eritaud, S. Schleimer]{Fran\c{c}ois Gu\'eritaud and Saul Schleimer}
%\date{February 2005}
\date{December 2006}

\begin{abstract}
Every cusped, finite--volume hyperbolic three--manifold has a ca\-nonical decomposition into ideal polyhedra. We study the canonical decomposition of the hyperbolic manifold obtained by filling some (but not all) of the cusps with solid tori: in a broad range of cases, generic in an appropriate sense, this decomposition can be predicted from that of the unfilled manifold.
\end{abstract}

\maketitle
%\tableofcontents

\section{Introduction}
\label{sec:intro}
%%%%%%%%%%%%%%%%%%%
%%               %%
%%   intro.tex   %%
%%               %%
%%%%%%%%%%%%%%%%%%%

Let $M$ be a complete cusped hyperbolic $3$--manifold of finite
volume, and endow the cusps $c_1,\dots,c_k$ of $M$ with disjoint
simple horoball neighborhoods $H_1,\dots,H_k$. The Ford--Voronoi
domain $\mathcal{F}\subset M$ consists of all points of $M$ 
%%% from which there is 
having a unique shortest path to the union of the $H_i$.  The
complement of $\mathcal{F}$ is a compact complex $C$ of totally
geodesic polygons. By definition, the canonical decomposition
$\mathcal{D}$ of $M$ with respect to the $H_i$ has one
$3$--dimensional cell (an ideal polyhedron) per vertex of $C$, one
face per edge of $C$, and one edge per (polygonal) face of $C$.  We
say that $\mathcal{D}$ is \emph{dual} to $C$. In
\cite{epstein-penner}, Epstein and Penner gave a precise description
of $\mathcal{D}$ in terms of convex hulls in Minkowski space
$\mathbb{R}^{3+1}$. Other names for $\mathcal{D}$ are the
geometrically canonical decomposition, or Delaunay (or Delone)
decomposition.

The combinatorics of $\mathcal{D}$, when the mutual volume ratios of
the $H_i$ are fixed, gives a complete topological invariant of the
manifold $M$, by Mostow's rigidity theorem. The canonical
decomposition $\mathcal{D}$ thus shows an interplay between
combinatorics on one hand, and hyperbolic geometry and
three--dimensional topology on the other. This motivates the study of
$\mathcal{D}$, and suggests that it is a difficult problem in
general. 
%%% Change in wording. 
Yet it seems that $\mathcal{D}$ can be determined in any
particular case by computation using, say, Jeffrey Weeks' computer
program SnapPea \cite{snappea}. General results are known only when
$M$ is restricted to belong to certain classes of manifolds:
punctured--torus bundles, two--bridge link complements, certain
arborescent link complements and related objects, or covers of any of
these spaces \cite{jorgensen, akiyoshi-canon, lackenby, aswy1, aswy2,
gufu, qf, these}. In fact, the combinatorics underlying all the above
examples are to a large extent the same. More examples, often using
symmetry, are compiled in \cite{sakuma-weeks}.

In the present paper, we will be interested in how the canonical
decomposition $\mathcal{D}$ changes when the last cusp $c_k$ (where
$k\geq 2$) undergoes a Dehn filling of slope $s$.  To this end, we
choose the reference horoball neighborhoods $\{H_i\}_{1\leq i <k}$ of
the remaining cusps after filling to have the same volumes as before
filling. Moreover, we make the running assumption that the horoball
neighborhood $H_k$, before filling, had much smaller volume than all
the other $H_i$ (by a result of Akiyoshi \cite{akiyoshi}, the
combinatorics of $\mathcal{D}$ can assume only a finite number of
different ``values'' as the volumes of the $\{H_i\}_{1\leq i \leq k}$
vary). Thurston showed that the metric of the Dehn filling converges,
in the sense of Gromov, to the metric of the unfilled manifold as the
filling slope goes to infinity 
%%% wording.
(choosing basepoints appropriately).  Accordingly, our philosophy will
be that as $c_k$ is filled, and thus replaced by a Margulis tube, only
the combinatorics inside or near the Margulis tube change, and only
in predictable fashion.  Cells away from the Margulis tube undergo
only a small geometric perturbation.

To ensure this, we will have to make the following ``genericity''
assumptions:

\begin{enumerate}[(I)]
\item The decomposition $\mathcal{D}$ (before filling) consists only
of ideal \emph{tetrahedra}; 
\item There exists a unique shortest path from $H_k$ to
$\bigcup_{i=1}^{k-1}H_i$ in $M$. 
%\item The filling slope in the cusp $c_k$ is sufficiently complicated
%(i.e. a finite number of ``forbidden'' slopes are avoided). 
\end{enumerate}

Of course, the above notion of genericity is problematic, since there
are only countably many complete finite--volume cusped hyperbolic
$3$--manifolds $M$ to choose from, and certainly, infinitely many of
them will be non--generic. Still, we have checked that 12 out of the
15 twice--cusped manifolds in the five--tetrahedron census of SnapPea
are generic.

\begin{theorem} \label{thm:main}
Under the genericity assumptions \emph{(I--II)} above, if the volume
of the cusp neighborhood $H_k$ is small enough, then the decomposition
$\mathcal{D}$ (before filling) contains exactly two ideal tetrahedra
$\Delta,\Delta'$ that have a vertex in the cusp $c_k$. Moreover,
$\Delta,\Delta'$ are isometric, each of $\Delta,\Delta'$ has exactly
one vertex in $c_k$, and $\partial (\Delta \cup \Delta')$ is a
once-punctured torus. For any sufficiently large filling slope $s$ in
the cusp $c_k$, the canonical decomposition $\mathcal{D}_s$ of the
manifold obtained by Dehn filling along $s$ is combinatorially of the
form
$$\mathcal{D}_s=\left ( \mathcal{D}\smallsetminus
\{\Delta,\Delta'\}\right ) \cup \mathcal{T}$$ where
$\mathcal{T}=\Delta_1\cup \dots \cup \Delta_N$ is a solid torus minus
one boundary point, and the combinatorial gluing of the $\Delta_i$ is
dictated by the continued fraction expansion of the slope $s$, with
respect to a certain basis of the first homology of the cusp $c_k$
depending only on $\mathcal{D}$.
\end{theorem}

Geometrically, the tetrahedra of $\mathcal{D}_s \smallsetminus
\mathcal{T}$ are small deformations of the tetrahedra of $\mathcal{D}
\smallsetminus (\Delta \cup \Delta')$. Section \ref{sec:farey} will
make explicit how the continued fraction expansion of $s$ dictates a
triangulation.  To predict $\mathcal{D}_s$ when genericity is not
satisfied, or even to estimate the number of slopes $s$ which fail to
be sufficiently large in the sense of Theorem \ref{thm:main} (their
number may not be universally bounded), remains very challenging.

We will prove Theorem \ref{thm:main} in Section
\ref{sec:voronoi}. Moreover, an analogous statement (Theorem
\ref{thm:main-mult}) will still hold when more than one cusp is
filled.  In Section \ref{sec:whitehead}, we will treat a real--life
family of examples by showing

\begin{theorem} \label{thm:whitehead}
If $M$ is a hyperbolic Dehn filling of one cusp of the Whitehead link
complement in $\mathbb{S}^3$, the canonical decomposition of $M$ is
dictated by the continued fraction expansion of the filling slope.
\end{theorem}

The Whitehead link complement actually violates both conditions of the
genericity assumption, but its symmetry compensates this
inconvenience. In fact, we will construct a certain triangulated solid
torus, also denoted $\mathcal{T}$, that serves as a proxy for the
Margulis (filling) tube itself: in the case of the Whitehead link
complement, it turns out that the filled manifold consists only of
$\mathcal{T}$ (with some exterior faces pairwise identified), i.e. no
combinatorics outside $\mathcal{T}$ need to be remembered from the
unfilled manifold. However, $\mathcal{T}$ can be slightly more
complicated than in Theorem \ref{thm:main} --- see Section
\ref{sec:whitehead} for details.

Historically, the first avatar of the triangulation $\mathcal{T}$ of
Theorem \ref{thm:main} seems to go back to \cite{jorgensen-cyclic}
where J{\o}rgensen briefly described the Ford--Voronoi domain of the
quotient of $\mathbb{H}^3$ by a loxodromy, with respect to an ideal
point. Full proofs of his results were since given, and the case of
non--ideal points included, by Drumm and Poritz in
\cite{drumm-poritz}.  Setting aside the case of non--ideal points, we
use \emph{angle structures} and ideal triangulations (combinatorially
dual to the Ford--Voronoi domain) to obtain new and quite different
proofs of their results.
%%% Before: Although it is arguable that our results on Dehn
%%% fillings might be reducible to their techniques, we use different
%%% methods (mainly the notion of \emph{angle structures}) and our basic
%%% objects are ideal triangulations, combinatorially dual to the
%%% Ford--Voronoi domain.  
Additionally, our paper provides the following improvements over the
existing literature:
%%% Before: Technical aspects aside, our paper thus provides
%%% the following (modest) factual improvements over the existing
%%% literature:
\begin{itemize}
\item Suppose that $\Gamma$ is a Kleinian group and $Z \subset \Gamma$ an
infinite cyclic subgroup resulting from a Dehn filling.  Then the
canonically triangulated solid torus corresponding to $Z$ is
incorporated into the canonical triangulation of
$\mathbb{H}^3/\Gamma$.
%%% New:
Under the genericity assumption this incorporation explains how, in
the program SnapPea, the picture of a triangulated cusp neighborhood
changes under Dehn filling.
\item In Section \ref{sec:elementary-extension}, we sketch an
extension to the case where $Z$ is virtually cyclic.
\item The convex hull in $\mathbb{H}^3$ of an ideal loxodromic orbit
always admits a canonical triangulation by the Epstein--Penner
construction (extended to the infinite--covolume case by Akiyoshi and
Sakuma \cite{aki-saku}). However, some of the outermost tetrahedra may
be \emph{timelike} or \emph{lightlike}, not \emph{spacelike}, in which
case they do not correspond to vertices of the Ford--Voronoi domain
(which indeed may have no vertices at all!). Although this case does
not arise in the context of Dehn fillings because the covolume stays
finite \cite{epstein-penner}, it is covered at no extra cost by our
methods, and apparently eludes those of \cite{drumm-poritz}.
\end{itemize}

The plan of the paper is as follows. In Section \ref{sec:rivin} we
recall the definition of the space $W$ of \emph{angle structures} on a
combinatorial ideal triangulation, and explain (following Rivin
\cite{rivin}) how to find the hyperbolic structure by maximizing a
volume functional $\mathcal{V}$ on $W$; an application is given for
solid tori. In Section \ref{sec:farey} we recall the combinatorics of
the Farey graph in $\mathbb{H}^2$ and use it to describe an ideal
triangulation of a solid torus $\mathcal{T}$. In Section
\ref{sec:voronoi}, using results from \cite{qf}, we check that the
decomposition of $\mathcal{T}$ is geometrically canonical, and
describe how to insert $\mathcal{T}$ as a proxy Margulis tube of a
filled manifold, under the ``genericity'' assumptions. In Section
\ref{sec:whitehead}, we adapt the method to treat all Dehn fillings on
one component of the Whitehead link complement.

\subsection*{Acknowledgements} We are very grateful to PCMI (Park
City), where this work originated during the summer of 2006. This
project would have been impossible without Jeff Weeks' program SnapPea
\cite{snappea}. 

This work is in the public domain.

\section{Angle structures and volume maximization} 
\label{sec:rivin}
%%%%%%%%%%%%%%%%%%%
%%               %%
%%   rivin.tex   %%
%%               %%
%%%%%%%%%%%%%%%%%%%

In Section \ref{sec:rivinstheorem} we give basic definitions and quote Theorem \ref{thm:rivin} (due to Rivin), the cornerstone of our method to find ideal triangulations. In Section \ref{sec:rigidity}, we parametrize the deformation space of certain hyperbolic solid tori; 
% (in particular this is a rigidity result)
the method, while not a direct application of Theorem \ref{thm:rivin}, follows from the same ideas and from the concept of ``spun'' triangulations \cite{thurstonnotes}.

\subsection{Rivin's theorem} \label{sec:rivinstheorem}

\begin{definition} A \emph{(combinatorial) ideal tetrahedron} is a space diffeomorphic to an ideal tetrahedron of hyperbolic space $\mathbb{H}^3$ (i.e. with vertices at infinity); the faces of such an ideal tetrahedron are called ideal triangles. \end{definition}

Consider an oriented combinatorial ideal tetrahedron $\Delta$, and copies $\Delta_1,\dots,\Delta_N$ of $\Delta$: the $\partial \Delta_i$ naturally receive consistent orientations. A \emph{gluing} of the $\Delta_i$ is an equivalence relation on $\bigsqcup_{i=1}^N \Delta_i$ generated by orientation--reversing identifications $\phi_{FG}:G\rightarrow F$ of pairs of faces $F\neq G$ of the $\Delta_i$, in such a way that 
\begin{itemize}
\item For each face $F$ of each $\Delta_i$, there is at most one face $G$ (resp. $G'$) of some $\Delta_j$ such that $\varphi_{FG}$ (resp. $\varphi_{G'F}$) is defined; moreover $G$ exists if and only if $G'$ exists and one then has $G=G'$ and $\varphi_{G'F}=\varphi_{FG}^{-1}$;
\item Whenever $\varphi:=\varphi_{F_1F_2}\circ \varphi_{F_2F_3}\circ\dots\circ \varphi_{F_{n-1}F_n}\circ\varphi_{F_nF_1}$ is well--defined on an edge $\epsilon$ of $\Delta_i$, then $\varphi$ is the identity of $\epsilon$.
\end{itemize}
The last condition is called the trivial holonomy condition.

Let $\sim$ be a gluing: then $M:=\left . \bigsqcup_{i=1}^N \Delta_i \right / \!\sim$ is a manifold (possibly with boundary). We say that the $\Delta_i$ endow $M$ with an ideal triangulation. The $6N$ edges of the $\Delta_i$ define edges in $M$, which we call boundary edges if they belong to $\partial M$, and interior edges otherwise.

Let us denote by $\epsilon^1_i, \dots,\epsilon^6_i$ the six edges of $\Delta_i$ (before gluing), and by $E$ the set of all $\epsilon^{\kappa}_i$ (so $|E|=6N$). We say that $\epsilon\in E$ is \emph{incident} to an edge $e$ of $M$ if $\epsilon$ projects to $e$ under the gluing ``$\sim$''. Fix a map $\alpha\colon \{ \text{boundary edges of }M \}\rightarrow \mathbb{R}_+^*$.

\begin{definition}
An \emph{angle structure} on $M$ with respect to $\alpha$ is a map $\theta:E\rightarrow \mathbb{R}_+^*$ such that 
\begin{itemize}
\item If the edges $\epsilon,\epsilon',\epsilon''$ of $\Delta_i$ share a vertex, then $\theta(\epsilon)+\theta(\epsilon')+\theta(\epsilon'')=\pi$;
\item If $\epsilon_1,\dots,\epsilon_n \in E$ is the full list of edges incident to an interior edge $e$ of $M$, then $\sum_{i=1}^n\theta(\epsilon_i)=2\pi$;
\item If $\epsilon_1,\dots,\epsilon_n \in E$ is the full list of edges incident to a boundary edge $e$ of $M$, then $\sum_{i=1}^n\theta(\epsilon_i)=\pi-\alpha(e)$.
\end{itemize}
\end{definition}
The $\theta(\epsilon)$, for $\epsilon \in E$, are called the dihedral angles of the $\Delta_i$. Given an angle structure, we can realize each $\Delta_i$ by an ideal hyperbolic tetrahedron $\delta_i$ of $\mathbb{H}^3$ with dihedral angles $\theta(\epsilon^1_i),\dots,\theta(\epsilon^6_i)$; however, when the face identifications $\varphi_{FG}$ are the corresponding hyperbolic isometries, the trivial holonomy condition may be violated. The following theorem tells us exactly for which angle structures this problem disappears.

\begin{theorem}[Rivin, \cite{rivin}] \label{thm:rivin}
Suppose the space $W$ of angle structures is non--empty. Then every critical point $\theta \in W$ of the volume functional 
%$$\begin{array}{rrcl}&W&\rightarrow&\mathbb{R}_+^* \\ \mathcal{V}:&\theta& \mapsto & -\frac{1}{2}\sum_{\epsilon \in E} \int_0^{\theta(\epsilon)} \log |2\sin u|\,du \end{array}$$ 
$$\mathcal{V}(\theta):= -\frac{1}{2}\sum_{\epsilon \in E} \int_0^{\theta(\epsilon)} \log |2\sin u|\,du ~>0$$ 
defines a complete hyperbolic metric with polyhedral boundary on $M$, with dihedral angle $\alpha(e)$ at each exterior edge $e$. Conversely, if $M$ admits such a complete hyperbolic metric in which the $\Delta_i$ are realized by totally geodesic ideal tetrahedra $\delta_i$ with disjoint interiors, then the dihedral angles of the $\delta_i$ define a critical point of $\mathcal{V}$.
\end{theorem}

Note that in an angle structure, the dihedral angles at opposite edges of any tetrahedron $\Delta_i$ are equal; if $\theta_1,\theta_2,\theta_3$ are the dihedral angles at the edges coming into one (and therefore any) vertex of $\Delta_i$, then $\mathcal{V}_0(\theta_1,\theta_2,\theta_3):=-\sum_{i=1}^3 \int_0^{\theta_i} \log |2\sin u|\,du$ is the volume of the ideal tetrahedron of $\mathbb{H}^3$ with those dihedral angles. This tetrahedron is unique up to isometry of $\mathbb{H}^3$.

\begin{fact} \label{fact:infinite-derivative}
The function $\mathcal{V}_0$ is convex on $\Theta:=\{(\theta_1,\theta_2,\theta_3)\in \mathbb{R}_+^3~|~\theta_1+\theta_2+\theta_3=\pi\}$, strictly convex on the interior of $\Theta$, and vanishes on $\partial \Theta$.
For any $x\in (0,\pi)$ and any $\omega \in \mathbb{R}$ one has 
$${\frac{d}{dt}}_{|t=0^+} \mathcal{V}_0(\pi-x-\omega t~,~x-(1-\omega)t~,~t)=+\infty.$$ 
%(when $x\in\{0,\pi\}$ however, the statement fails for all values of $\omega$).
This expresses the fact that if exactly one angle of an ideal tetrahedron $\Delta$ is $0$, increasing that angle to $\varepsilon<<1$ yields a volume increase much greater than $\varepsilon$; note that the same statement is false when \emph{two} angles of $\Delta$ are $0$.
\end{fact} 

Fact \ref{fact:infinite-derivative} implies that the volume functional $\mathcal{V}:W\rightarrow \mathbb{R}$ of Theorem \ref{thm:rivin} is concave and positive, and extends continuously to a concave function on the (compact) closure $\overline{W}$ of $W$. It moreover implies

\begin{proposition}[Rivin, \cite{rivin}] \label{prop:degeneracies} Suppose $W\neq \emptyset$ and let $\theta_0 \in \overline{W}$ be a point where $\mathcal{V}$ reaches its maximum. Either
\begin{itemize}
\item $\theta_0$ belongs to $W$, i.e. $\theta_0(E)\subset \mathbb{R}_+^*$, in which case $\theta_0$ is a (necessarily unique) critical point for $\mathcal{V}$ in $W$; or
\item there exists a non--empty list of tetrahedra $\Delta_{i_1},\dots,\Delta_{i_s}$ that have an edge $\epsilon$ such that $\theta_0(\epsilon)=0$: then, each $\Delta_{i_k}$ also has an edge $\epsilon'$ such that $\theta_0(\epsilon')=\pi$.
\end{itemize}
\end{proposition}

%%%%%%%%%%%%%%%%%%%%%%%%%%%%%%%%%%%%

\subsection{Rigidity of solid tori} \label{sec:rigidity}

In this section we prove a rigidity result for hyperbolic polyhedral solid tori with given dihedral angles (and one ideal vertex). The method is a special case of a generalization of Theorem \ref{thm:rivin} to spun triangulations.

Consider a once--punctured torus $\tau$ with three ideal edges $e,e',e''$ running from the puncture to itself: these edges divide $\tau$ into two ideal triangles. Let $\gamma$ be a non--oriented free homotopy class of simple closed curves in $\tau$, and let $n,n',n'' \in \mathbb{N}$ be the minimal intersection numbers of $\gamma$ with $e,e',e''$ respectively. It is well--known that the triple $(n,n',n'')$ determines the class $\gamma$, and that the largest among $n,n',n''$ is the sum of the other two terms.

Let $a,b,c \in [0,\pi)$ be such that $a+b+c=\pi$. 
We aim to construct a punctured solid torus $X$ (namely a solid torus minus one point of its boundary) with the following properties: the punctured torus $\partial X$ has three ideal edges with exterior dihedral angles $a,b,c$, and there exist coprime positive integers $n_a, n_c$ such that the meridian of $X$ intersects these three edges minimally in $n_a,\,n_a+n_c,\,n_c$ points respectively. We write $n_b:=n_a+n_c$. 

\begin{proposition}
A hyperbolic solid torus $X$ as above exists if and only if $a\,n_a+b\,n_b+c\, n_c>2\pi$. This solid torus is then unique up to isometry. \label{prop:rigidity}
\end{proposition}
\begin{remark}
The left member of the inequality is the sum of exterior dihedral angles met by a meridian in $\partial X$: the inequality can thus be seen as a sort of Gauss-Bonnet condition for the compression disk of the solid torus $X$ (see \cite{fugu} for a more general construction). In Section \ref{sec:farey}, we will check that the same condition is enough for a certain (non--spun) ideal triangulation of $X$ to have angle structures (with respect to $a,b,c$), and indeed to be geometrically realized.
\end{remark}
\begin{proof}
If $X$ exists, we can consider its universal cover $U$ which is a complete hyperbolic manifold with locally convex boundary and is thus, by a standard argument, naturally embedded in $\mathbb{H}^3$. This space $U$ is the convex hull of the orbit of an ideal point of $\partial_{\infty}\mathbb{H}^3\simeq \mathbb{S}^2$ under a certain loxodromic $\varphi$ (corresponding to the core curve of $X$). We can stellate $U$ with respect to the attractive fixed point of $\varphi$: this yields a $\varphi$--invariant decomposition of $U$ (minus the axis of $\varphi$) into tetrahedra, hence, quotienting out by $\varphi$, a decomposition of the solid torus $X$ (minus the core axis) into two ideal tetrahedra $\Delta,\Delta'$. Note that this decomposition has only one interior edge $L$, originating at the puncture of $\partial X$ and spinning towards the core of $X$. Thus, constructing $X$ in general amounts to finding positive dihedral angles for $\Delta,\Delta'$ such that 
\begin{enumerate}[(i)]
\item the holonomy around $L$ is the identity of $\mathbb{H}^3$, i.e. the complex angles around $L$ sum to $2\pi$;
% COMPLEX ANGLES : NOT DEFINED....
\item the boundary of $\Delta \cup \Delta'$ has interior dihedral angles $\pi-a,\, \pi-b,\,\pi-c$;
\item the holonomy around the core curve of $X$ is also the identity of $\mathbb{H}^3$.
\end{enumerate}

(One may refer e.g. to Definition 6.3 of \cite{gufu} for a precise definition of holonomy.)
Condition (i) above is automatically satisfied because each dihedral angle of $\Delta$ and $\Delta'$ is incident to $L$ exactly once. To study Condition (ii), let us fix some notation: let $ABC$ and $ACD$ be two counterclockwise oriented triangles in $\mathbb{C}\subset \mathbb{P}^1\mathbb{C} \simeq \partial_{\infty}\mathbb{H}^3$; we identify $\Delta$ with the tetrahedron $\infty ABC$ and $\Delta'$ with $\infty ACD$, gluing the ideal triangles $\infty AB$ and $\infty DC$ (resp. $\infty AD$ and $\infty BC$) together. The interior angles at $A,B,C$ of $\Delta$ are noted $\delta_a, \delta_b, \delta_c$ respectively. The interior angles at $A,C,D$ of $\Delta'$ are noted $\delta'_c, \delta'_a, \delta'_b$ respectively (see Figure \ref{fig:holonomies}). 

\begin{figure}[h!] \centering
\psfrag{A}{$A$}
\psfrag{B}{$B$}
\psfrag{C}{$C$}
\psfrag{D}{$D$}
\psfrag{da}{$\delta_a$}
\psfrag{db}{$\delta_b$}
\psfrag{dc}{$\delta_c$}
\psfrag{dpa}{$\delta'_a$}
\psfrag{dpb}{$\delta'_b$}
\psfrag{dpc}{$\delta'_c$}
\psfrag{la}{$\lambda_a$}
\psfrag{lc}{$\lambda_c$}
\includegraphics{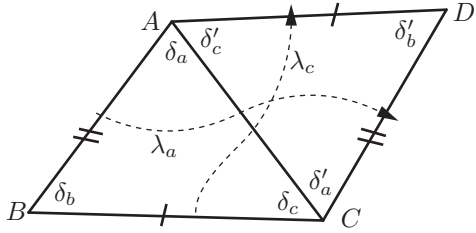}
\caption{The cusp shapes of $\Delta$ and $\Delta'$.} \label{fig:holonomies} 
% INCLUDE A PICTURE FROM THE OTHER CUSP (A HEXAGON)
\end{figure}

Condition (ii) can then be written 
$$\delta_a+\delta'_a=\pi-a~;~\delta_b+\delta'_b=\pi-b~;~\delta_c+\delta'_c=\pi-c.$$ 
This implies 
\begin{equation} \label{eq:Delta-angles} \left \{ \begin{array}{rcl}(\delta_a,\delta_b, \delta_c)&=& \left ( \frac{\pi-a}{2}+\alpha~,~\frac{\pi-b}{2}+\beta~,~ \frac{\pi-c}{2}+\gamma \right ) \\ (\delta'_a,\delta'_b,\delta'_c) &=& \left ( \frac{\pi-a}{2}-\alpha~,~\frac{\pi-b}{2}-\beta~,~ \frac{\pi-c}{2}-\gamma \right ) \end{array} \right . \end{equation} 
where 
\begin{equation} \label{eq:hexagon}
\textstyle{|\alpha|<\frac{\pi-a}{2}~,~|\beta|<\frac{\pi-b}{2}~,~|\gamma|<\frac{\pi-c}{2}}~,~\text{ and } \alpha+\beta+\gamma=0. \end{equation}
The space of solutions $(\alpha,\beta,\gamma)$ to (\ref{eq:hexagon}) is the interior of a centrally symmetric affine hexagon $P$ whose edges are given by 
\begin{equation} \label{eq:6edges}
\textstyle{\alpha=\frac{\pi-a}{2}~,~\beta=-\frac{\pi-b}{2}~,~\gamma=\frac{\pi-c}{2}~,~\alpha=-\frac{\pi-a}{2}~,~\beta=\frac{\pi-b}{2}~,~\gamma=-\frac{\pi-c}{2}}
\end{equation}
in that order. (It is easy to check that these edges are all non--empty segments if $a,b,c>0$, and that e.g. the first and fourth edges are reduced to points if and only if $a=0$.) See Figure \ref{fig:davidstar}.

\begin{figure}[h!] \centering
\psfrag{ap}{$\alpha=\frac{\pi-a}{2}$}
\psfrag{am}{$\alpha=-\frac{\pi-a}{2}$}
\psfrag{bp}{$\beta=\frac{\pi-b}{2}$}
\psfrag{bm}{$\beta=-\frac{\pi-b}{2}$}
\psfrag{cp}{$\gamma=\frac{\pi-c}{2}$}
\psfrag{cm}{$\gamma=-\frac{\pi-c}{2}$}
\psfrag{P}{$P$}
\psfrag{S}{$S$}
\includegraphics{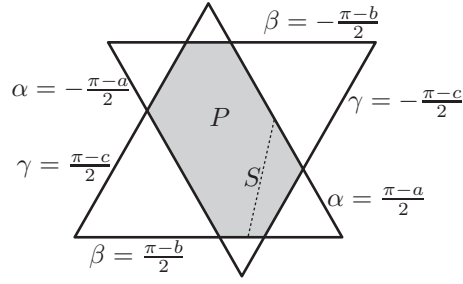}
\caption{The hexagon $P$ of solutions $(\alpha,\beta,\gamma)$ to (\ref{eq:hexagon}), and the segment $S$ of spun angle structures. At most one pair of opposite sides of $P$ can be reduced to points, because $a,b,c<\pi$.} \label{fig:davidstar} \end{figure}

Condition (iii) has two components: first, an angular component (affine in terms of the dihedral angles of $\Delta,\Delta'$) which will narrow down the space of solutions $(\alpha,\beta,\gamma)$ to the intersection of the interior of $P$ with a certain line. This intersection will be non--empty (namely, an open segment $S$) exactly when the inequality of Proposition \ref{prop:rigidity} is satisfied. Second, a scaling component which we will solve by seeking a critical point of a volume functional on $S$.

{\bf Angular component.} Following the notation above (and Figure \ref{fig:holonomies}), we refer to the three exterior edges of $X$ as $AB, BC, CA$: the corresponding exterior dihedral angles are $c,a,b$ respectively. 
%Recall that the meridian $\mu$ of $X$ intersects $AB, BC, CA$ respectively $n_c$, $n_a$, and $n_c+n_a$ times.
The angular holonomy map is a group homomorphism $h:H_1(\partial X, \mathbb{Z}) \rightarrow \mathbb{R}$. Let the oriented closed curve $\lambda_a$ (resp. $\lambda_c$) be a boundary component of a regular neighborhood of the oriented edge $\overrightarrow{BC}$ (resp. $\overrightarrow{BA}$), as in Figure \ref{fig:holonomies}. By the conventions (\ref{eq:hexagon}) above, $h([\lambda_a])=\delta_a-\delta'_a=2\alpha$ and $h([\lambda_c])=-\delta_c+\delta'_c=-2\gamma$. The meridian $\mu$ of $X$ is homotopic to $n_c [\lambda_a] + n_a [\lambda_c]$ (where $n_a,n_c>0$), because this class intersects $n_a$ times the edge $BC$, $n_b=n_a+n_c$ times the edge $AC$, and $n_c$ times the edge $BA$. Hence, $$h([\mu])=2(n_c\alpha-n_a\gamma).$$ Using (\ref{eq:6edges}), and considering the appropriate vertex of the space $P$ of angle structures, the largest (resp. smallest) possible value of $h([\mu])$ on the closure of $P$ is therefore $2\left ( n_c\frac{\pi-a}{2}+n_a\frac{\pi-c}{2}\right)=a n_a+b n_b+c n_c$ (resp. the negative of that number), where we used $a+b+c=\pi$ and $n_b=n_a+n_c$. We conclude that $h([\mu])=2\pi$ is satisfiable on the interior of $P$ if and only if $a\,n_a+b\,n_b+c\,n_c>2\pi$, as wished.

{\bf Scaling component.} The scaling holonomy map is a group homomorphism $\eta:H_1(\partial X,\mathbb{Z})\rightarrow \mathbb{R}_+^*$. The sine formula for triangles yields $\eta([\lambda_a])=\frac{\sin \delta_b}{\sin \delta_c}\,\frac{\sin \delta'_c}{\sin \delta'_b}$ and $\eta([\lambda_c])=\frac{\sin \delta_b}{\sin \delta_a}\,\frac{\sin \delta'_a}{\sin \delta'_b}$, hence $$\eta([\mu])=\eta([\lambda_a])^{n_c}\eta([\lambda_c])^{n_a}=\left ( \frac{\sin \delta_b}{\sin \delta'_b}\right )^{n_a+n_c} 
\left ( \frac{\sin \delta_c}{\sin \delta'_c}\right )^{-n_c}
\left ( \frac{\sin \delta_a}{\sin \delta'_a}\right )^{-n_a}.$$

On the other hand, let $S$ be the open segment defined by the intersection of the interior of $P$ with the condition $h([\mu])=2\pi$, i.e. $n_c\alpha-n_a\gamma=\pi$. The tangent space of $S$ is generated by the vector $(\dot{\alpha},\dot{\beta},\dot{\gamma})=(n_a,-n_a-n_c,n_c)$. Let $\Lambda$ be the Lobachevski function defined by $\Lambda(x)=-\int_{0}^{x}\log |2 \sin t|\, dt$. The volume functional is by definition 
$$\begin{array}{rrcl}&S& \longrightarrow & \mathbb{R}^+ \\ 
\mathcal{V}: & (\alpha,\beta,\gamma) & \longmapsto & 
\mathcal{V}_0(\delta_a,\delta_b,\delta_c) + \mathcal{V}_0(\delta'_a,\delta'_b,\delta'_c) \end{array}$$
where $\mathcal{V}_0(x,y,z)=\Lambda(x)+\Lambda(y)+\Lambda(z)$ is the volume of one ideal tetrahedron, and $\delta_a,\dots,\delta'_c$ are given by (\ref{eq:Delta-angles}). By Fact \ref{fact:infinite-derivative}, $\mathcal{V}$ is strictly concave on the segment $S$ and achieves its maximum in $S$ (indeed, the endpoints of $S$ belong to the perimeter of the hexagon $P$, but at any point of $\partial P$, at least one of the tetrahedra $\Delta,\Delta'$ has exactly one angle whose value is $0$: therefore, $\mathcal{V}$ has unbounded derivative near each endpoint of $S$). As a result, $\mathcal{V}$ has a unique (critical) maximum in the open segment $S$.

At that critical point, since $(\dot{\alpha},\dot{\beta},\dot{\gamma})=(n_a~,~-n_a-n_c~,~n_c)$, we have 
\begin{eqnarray*}
0=d\mathcal{V}(\dot{\alpha},\dot{\beta},\dot{\gamma}) &=&
\dot{\alpha}\Lambda'(\delta_a)+\dot{\beta}\Lambda'(\delta_b)+\dot{\gamma}\Lambda'(\delta_c)-\dot{\alpha}\Lambda'(\delta'_a)-\dot{\beta}\Lambda'(\delta'_b)-\dot{\gamma}\Lambda'(\delta'_c) \\ &=&
-\dot{\alpha}\log |2 \sin \delta_a|
-\dot{\beta}\log |2 \sin \delta_b|
-\dot{\gamma}\log |2 \sin \delta_c| \\ && 
+\dot{\alpha} \log |2 \sin \delta'_a| 
+\dot{\beta} \log |2 \sin \delta'_b|
+\dot{\gamma} \log |2 \sin \delta'_c| \\ &=& \log \left [ 
\left ( \frac{\sin \delta_a}{\sin \delta'_a} \right )^{-n_a}\, 
\left ( \frac{\sin \delta_b}{\sin \delta'_b} \right )^{n_a+n_c} \,
\left ( \frac{\sin \delta_c}{\sin \delta'_c} \right )^{-n_c} \right ] 
=\log \eta([\mu]). 
\end{eqnarray*}

At the critical point of $\mathcal{V}$ in $S$, we therefore have the following values for the holonomy maps: $h([\mu])=2\pi$ (rotational component) and $\eta([\mu])=1$ (scaling component). This precisely means that the metric completion of $\Delta\cup \Delta'$ is the solid torus $X$ endowed with a spun triangulation of two tetrahedra whose tips spin around the core curve. Moreover, since the critical point of $\mathcal{V}$ in $S$ is unique, we have in fact proved that $X$ is unique up to isometry. \end{proof}

\section{Farey combinatorics in solid tori} 
\label{sec:farey}
%%%%%%%%%%%%%%%%%%%
%%               %%
%%   farey.tex   %%
%%               %%
%%%%%%%%%%%%%%%%%%%

Let $X$ be a compact solid torus, minus one point of its boundary; call this removed point the \emph{puncture}.

%\begin{definition} A \emph{(combinatorial) ideal tetrahedron} is a space diffeomorphic to an ideal tetrahedron of hyperbolic space $\mathbb{H}^3$ (i.e. with vertices at infinity); the faces of such an ideal tetrahedron are called ideal triangles. \end{definition}

In this section we will first describe a certain combinatorial decomposition $\mathcal{D}$ of $X$ into ideal tetrahedra, relative to a given ideal triangulation of $\partial X$ (into two ideal triangles). This decomposition has previously been described, and studied in great detail, by Jaco and Rubinstein \cite{jaco-rubinstein}. We will then go on to find a geometric realization of $\mathcal{D}$, using the ideas of Section \ref{sec:rivin}.

\subsection{The Farey graph}

Identify the boundary at infinity of the hyperbolic plane $\mathbb{H}^2$ to the circle $\mathbb{P}^1\mathbb{R}$, endowed with the action of $PSL_2(\mathbb{Z})$. We assume that $0,1,\infty$ lie counterclockwise in that order on $\partial_{\infty}\mathbb{H}^2\simeq \mathbb{P}^1\mathbb{R}$. Consider the subset $\mathbb{P}^1\mathbb{Q}$ of $\mathbb{P}^1\mathbb{R}$. We measure the ``proximity'' of two elements $q=\frac{y}{x}$ and $q'=\frac{y'}{x'}$ of $\mathbb{P}^1\mathbb{Q}$ (given as ratios of coprime integers) by computing their wedge 
\begin{equation} \label{eq:wedge}
q\wedge q' :=\left | \left | \begin{array}{cc} y & y' \\ x & x' \end{array} \right | \right | \: \in \mathbb{N} \hspace{20pt} \text{(absolute value of the determinant).}\end{equation}
If we draw a straight line in $\mathbb{H}^2$ from $q$ to $q'$ each time $q\wedge q'=1$, we obtain the \emph{Farey triangulation} of $\mathbb{H}^2$. Alternatively, this triangulation can be defined by reflecting the ideal triangle $1\infty 0$ in its sides \emph{ad infinitum}.

Fix an identification (homeomorphism) between the punctured torus $\partial X$ and $\mathbb{T}:=(\mathbb{R}^2\smallsetminus \mathbb{Z}^2)/\mathbb{Z}^2$. We assume that the canonical orientation of $\mathbb{T}$ (induced by $\mathbb{R}^2$), followed by the outward--pointing normal of $\partial X$, coincides with the positive orientation on $X$.  The segment from $(0,0)$ to $(x,y)$ in $\mathbb{R}^2$ (where $x,y$ are coprime integers) projects to a properly embedded line $\gamma$ in $\partial X$: we say that $\frac{y}{x}\in \mathbb{P}^1\mathbb{Q}$ is the \emph{slope} of $\gamma$. An edge $E$ of the Farey triangulation (or: a Farey edge) corresponds to a pair of disjoint lines $\gamma,\gamma'$ in $\partial X$, whose slopes are the two ends of $E$ in $\mathbb{P}^1\mathbb{Q}$, and whose complement in $\partial X$ is an ideal quadrilateral. Similarly, Farey triangles (such as $1\infty 0$), having three vertices in $\mathbb{P}^1\mathbb{Q}$, correspond to triples of disjoint lines $\gamma,\gamma',\gamma''$ in $\partial X$ which define a decomposition of $\partial X$ into two ideal triangles. Finally, note that we can also associate a slope in $\mathbb{P}^1\mathbb{Q}$ to the meridinal closed curve $\mu$ of the solid torus $X$: namely, the slope of the unique properly embedded line $\mu'$ which (possibly after isotopy) does not intersect $\mu$.

Let $pqr$ be a Farey triangle, and suppose $m\in \mathbb{P}^1\mathbb{Q} \smallsetminus \{p,q,r\}$ is the slope of the meridian of $X$. By convention, we will suppose that the Farey edge $pq$ separates $r$ from $m$, and that $pqm$ is not a Farey triangle (so $m$ is ``far enough'' from the triangle $pqr$). Endow the punctured torus $\partial X$ with the ideal triangulation associated to $pqr$ (which we call the $pqr$--triangulation). In Section \ref{sec:triangulation}, we will be preoccupied with decomposing $X$ into ideal tetrahedra with faces (ideal triangles) glued in pairs, in such a way that exactly two ideal triangles remain free, and give the $pqr$-triangulation of $\partial X$.

\subsection{An ideal triangulation of the solid torus} \label{sec:triangulation}

The idea is to follow a path $\ell$ in the Farey triangulation, transverse to the Farey edges, from the ideal vertex $r$ to the ideal vertex $m$. We assume that the path $\ell$ crosses each Farey triangle at most once, i.e. never backtracks. The sequence of Farey triangles that $\ell$ encounters is then completely determined (so we can take $\ell$ to be e.g. a geodesic ray): these triangles are $$(T_0, T_1, \dots, T_N)=(pqr, pqr',\dots, mst)$$ where $s,t$ belong to $\mathbb{P}^1\mathbb{Q}$ and the symmetry of axis $pq$ takes $r$ to $r'$. Note that by assumption, $N\geq 2$.

For each $0<i<N$, we can then consider a properly embedded punctured torus $\tau_i\subset X$ isotopic to $\partial X$ (properness here means that by intersectiong $\tau_i$ with a basis of neighborhoods of the puncture of $X$, we get a basis of neighborhoods of the puncture of $\tau_i$). We can assume that the $\tau_i$ are disjoint and that $\tau_i$ separates $\partial X$ from $\tau_{i+1}$ (i.e. $\tau_{i+1}$ lies in $X$ ``inward'' from $\tau_i$). Endow $\tau_i$ with the triangulation associated to the Farey triangle $T_i$ --- for that purpose we also rechristen $\partial X$ as $\tau_0$. Note that two consecutive punctured tori $\tau_{i-1}, \tau_i$ always have two edge slopes in common (these slopes are the ends of the Farey edge $T_{i-1}\cap T_i$). Thus, we can isotope $\tau_1$ until its edges of slopes $p,q$ coincide with those of $\tau_0=\partial X$; then isotope $\tau_2$ until two of its edges coincide with the edges of similar slopes in $\tau_1$; then isotope $\tau_3$ until it intersects $\tau_2$ along two edges, etc.

At the end of this process, the space comprised between $\tau_{i-1}$ and $\tau_i$, for each $0<i<N$, is a (combinatorial) ideal tetrahedron $\Delta_{i}$ with four of its edges identified in opposite pairs. These tetrahedra $\Delta_{i}$, with the combinatorial gluing that arises from the construction above, are those of our decomposition $\mathcal{D}$ of $X$. (Since $N\geq 2$, there is at least one tetrahedron $\Delta_i$. Our ``half--shift'' convention $\partial \Delta_i=\tau_{i-1}\cup\tau_i$, or equivalently $\tau_i=\Delta_i\cap\Delta_{i+1}$, is arbitrary). In order to homotopically ``kill'' the meridian of the solid torus $X$, it only remains to describe the gluing of the last surface $\tau_{N-1}$ to itself.

\begin{figure}[h!] \centering
\psfrag{a}{$A$}
\psfrag{b}{$B$}
\psfrag{c}{$C$}
\psfrag{m}{$\mu$}
\includegraphics{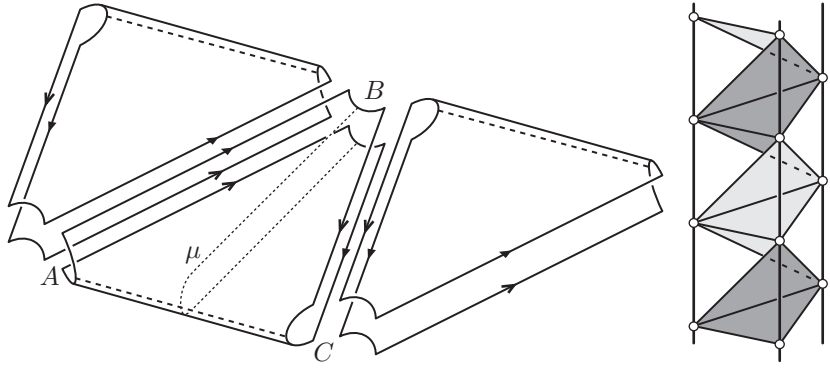}
\caption{Left: a punctured torus (shown are 3 folded copies of a fundamental domain; arrows are identified) wraps around an ideal M\"obius band. The meridian line $\mu$, of slope $m$, becomes homotopically trivial. The dotted folding line $AC$ has slope $m'$. Right: part of the universal cover of the same M\"obius band (shaded) and the tetrahedron $\Delta_{N-1}$ glued to it.} \label{fig:wrap} \end{figure}

If $T_N=mst$ is the last Farey triangle, let $T_{N-1}=m'st$ be the next-to-last, associated to the surface $\tau_{N-1}$. We fold $\tau_{N-1}$ along its edge of slope $m'$, gluing the two adjacent faces (ideal triangles) $F',F''$ to one another to obtain a single ideal triangle $F$. Intrinsically, $F$ is an ideal M\"obius band, i.e. a compact M\"obius band minus one point of its boundary. Indeed, from an (ideal) triangle $ABC$, one can construct an (ideal) M\"obius band $F$ with boundary $AC$, by gluing the oriented edge $AB$ to $BC$: the (punctured) torus $\tau_{N-1}=F'\cup F''$ then just wraps around this (ideal) M\"obius band $F$, like the boundary of a regular neighborhood of an embedding of $F$ in $\mathbb{R}^3$. See Figure \ref{fig:wrap}.

\subsection{Angle structures} \label{sec:angle-structures}

We proceed to describe positive angle structures for the tetrahedra $\Delta_{i}$, where $1\leq i \leq N-1$ (the argument is reminiscent of \cite{gufu} and \cite{qf}, although the solution space will look quite different).
%To simplify notation, we choose (somewhat arbitrarily) to write $\Delta_i:=\Delta_{i}$. 
%%This choice is somewhat arbitrary, but only the numbering of the $\Delta_i$ (as opposed to the surfaces $\tau_i$) will matter from here on. 
More precisely, consider reals 
\begin{equation} \label{eq:thetas}
\theta_p, \theta_q,\theta_r~\text{ such that } 
\left \{ \begin{array}{rcl}
\theta_p+\theta_q+\theta_r &=& \pi~; \\
\theta_p~,~\theta_q &\geq & 0~; \\
\pi~>~\theta_r &>& 0~. \end{array} \right .
\end{equation}
We will look for angle structures on the $\Delta_{i}$ such that the interior dihedral angles of $X$ at the edges of slope $p,q,r$ in $\partial X$ are $\pi-\theta_p, \pi-\theta_q, \pi-\theta_r$ respectively. Note that we do not allow $\theta_r$ to vanish: indeed, $\pi-\theta_r$ will be a dihedral angle of the first tetrahedron $\Delta_{1}$. (If the solid torus $X$ admits a geometric realization in which $\theta_r=0$, we can always remove this flat tetrahedron $\Delta_{1}$ and see $\partial X$ as being endowed with the $pqr'$--triangulation, where $r'\in\mathbb{P}^1\mathbb{Q}$ is the symmetric of $r$ with respect to the Farey edge $pq$.)

\begin{proposition} \label{prop:pqr}
An angle structure satisfying (\ref{eq:thetas}), also called a $(\theta_p,\theta_q, \theta_r)$--angle structure, exists if and only if $(m\wedge p) \theta_p + (m\wedge q)\theta_q + (m\wedge r) \theta_r >2\pi$.
\end{proposition}

%\begin{remark} The left member of the inequality is the sum of exterior dihedral angles met by a simple closed curve of slope $m$ in $\partial X$: the inequality can thus be seen as a sort of Gauss-Bonnet inequality for the compression disk of the solid torus $X$ (see \cite{fugu} for a more general construction). In Section \ref{sec:maximize}, we will check that the same condition is enough for an actual geometric triangulation with exterior angles $\theta_p, \theta_q, \theta_r$ to exist. \end{remark}

\begin{figure}[h!] \centering
\psfrag{po}{$(\pi,0,0)$}
\psfrag{opo}{$(0,\pi,0)$}
\psfrag{op}{$(0,0,\pi)$}
\psfrag{9}{$9$}
\psfrag{10}{$10$}
\psfrag{11}{$11$}
\psfrag{12}{$12$}
\psfrag{13}{$13$}
\psfrag{14}{$14$}
\psfrag{15}{$15$}
\psfrag{16}{$16$}
\psfrag{17}{$17$}
\psfrag{18}{$18$}
\psfrag{19}{$19$ or more}
\includegraphics{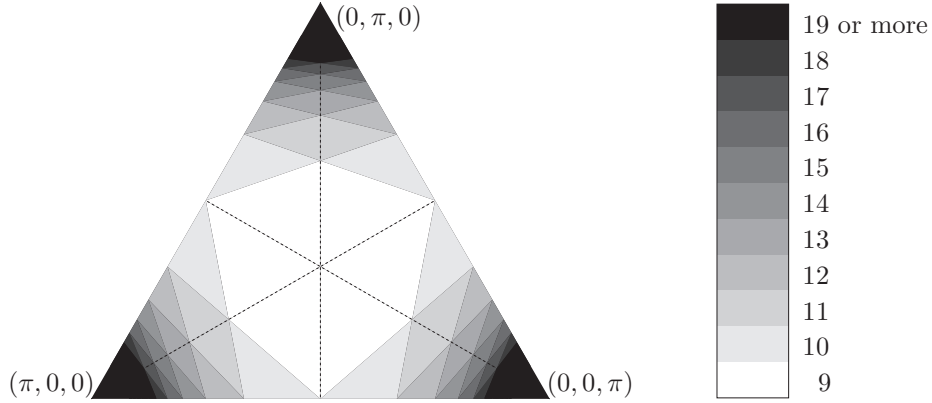}
\caption{Parameter space for the triple $(\theta_p, \theta_q, \theta_r)$, and numbers of ``forbidden'' slopes $m$ (the brighter, the fewer).} \label{fig:zebre} \end{figure}

\begin{remark} \label{rem:zebre}
It is easy to check that $m\wedge r = (m\wedge p)+(m\wedge q)$ --- e.g. by reducing to the case $(p,q)=(0,\infty)$ and using the $PSL_2(\mathbb{Z})$--invariance of the $\wedge$--notation. Thus, by (\ref{eq:thetas}), the inequality of Proposition \ref{prop:pqr} is automatically true unless $\text{min}\,\{m\wedge p , m\wedge q \}=1$. For example, if $(m\wedge p, m\wedge q)=(1,1)$, the condition is always false (recall we required that $pqm$ not be a Farey triangle); if $(m\wedge p, m\wedge q)=(2,1)$, it amounts to $\theta_r>\theta_q$. The equilateral triangle in Figure \ref{fig:zebre} shows the full parameter space for the triple $(\theta_p, \theta_q, \theta_r)$: shades indicate how many slopes $m$ fail to satisfy the condition of Proposition \ref{prop:pqr}, where we allow $m$ to range over all of $\mathbb{P}^1\mathbb{Q}$ rather than just over the arc $\overset{\frown}{pq}$ (when $m$ belongs to one of the arcs $\overset{\frown}{qr},\overset{\frown}{rp}$, we construct the same ideal triangulations, up to a permutation of $p,q,r$).
\end{remark}

\begin{proof} (Prop. \ref{prop:pqr}).
The tetrahedra $\Delta_{i}$ are naturally associated to the Farey edges $e_{i}=T_{i-1}\cap T_i$ that the path $\ell$ crosses. Orient $\ell$ from $T_0$ to $T_N$. If $e_{i}$ and $e_{i+1}$ share their Right (resp. Left) end with respect to the orientation of $\ell$, we say that $\ell$ \emph{makes a Right} (resp. \emph{a Left}) between $e_{i}$ and $e_{i+1}$ (or: at $T_i$). Thus, $\ell$ defines a word $\Omega=RLL...R$ of length $N-1$ in the letters $R,L$: for each $i\in \{1,2,\dots, N-1\}$ there is a tetrahedron $\Delta_{i}$ and a letter $\Omega_i\in\{R,L\}$. If $(p,q,r)=(0,\infty,-1)$, then the lengths of the syllables $R^n$ and $L^n$ of $\Omega$ are exactly the integers in the continued fraction expansion of the rational $m$, as referred to in Theorem \ref{thm:main}.

Note that no letter $R$ or $L$ is associated to the very first Farey triangle $T_0=pqr$, because the line $\ell$ does not ``enter'' $T_0$ through $pr$ rather than through $qr$. We nevertheless decide to place an extra letter $\Omega_0\in\{R,L\}$ in front of the word $\Omega$, so that $\Omega$ becomes of length $N$ and starts with either $RR$ or $LL$. This convention is totally artificial (the other choice would be equally good), but making a choice here will allow us to streamline the notation in our argument. Up to switching $p$ and $q$, we can now assume that $\ell$ enters the Farey triangle $T_0$ through the edge $pr$, and leaves through $pq$. See Figure \ref{fig:fareyline}

\begin{figure}[h!] \centering
\psfrag{p}{$p$}
\psfrag{q}{$q$}
\psfrag{r}{$r$}
\psfrag{rr}{$r'$}
\psfrag{m}{$m$}
\psfrag{s}{$s$}
\psfrag{t}{$t$}
\psfrag{mm}{$m'$}
\psfrag{T0}{$T_0$}
\psfrag{T1}{$T_1$}
\psfrag{TN}{$T_N$}
\psfrag{TNm}{$T_{N-1}$}
\psfrag{R}{$R$}
\psfrag{L}{$L$}
\psfrag{el}{$\ell$}
\includegraphics{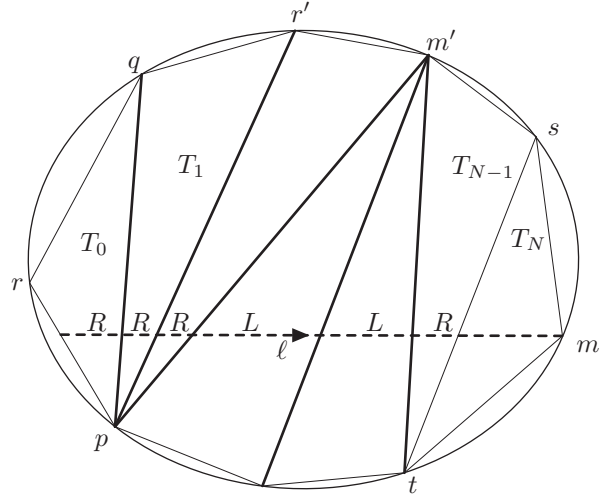}
\caption{The Farey graph. The $5$ thick lines $T_{i-1}\cap T_i$ (where $1\leq i \leq 5$) correspond to the tetrahedra $\Delta_{i}$.} \label{fig:fareyline} 
\end{figure}

\begin{definition}
If $\Omega_{i-1}\neq \Omega_i$, we say that $\Delta_{i}$ is a \emph{hinge} tetrahedron. Otherwise, we call $\Delta_{i}$ \emph{non--hinge}. For example, following our convention, $\Delta_{1}$ is non--hinge.
\end{definition}

To compute angle structures, it will be useful to describe the cusp triangulation associated to the ideal triangulation $\{\Delta_{i}\}_{1\leq i \leq N-1}$ of $X$. Since each pleated punctured torus $\tau_i$ has one ideal vertex and three edges, each with two ends, the link of the ideal vertex of $\tau_i$ is a hexagon $H_i$ (the pleating angles of $\tau_i$ are the exterior angles of $H_i$). We are going to define the dihedral angles of the ideal tetrahedra $\Delta_{i}$ in terms of the pleating angles of the $\tau_i$. Note that the hexagon $H_i$ has a central symmetry induced by the hyperelliptic involution of the punctured torus $\tau_i$ (rotation of $180^{\circ}$ around the puncture, which exchanges the ends of each edge of $\tau_i$).

Let $\xi\eta\zeta=T_{i-1}$ and $\xi\eta\zeta'=T_i$ be two consecutive Farey triangles, so that the Farey vertex $\xi$ (resp. $\eta$) lies to the right (resp. left) of the oriented axis $\ell$. The tetrahedron $\Delta_{i}$ has: 
\begin{itemize}
\item two opposite edges carrying the same dihedral angle $x_{i}$ and identified to just one edge, of slope $\xi$, in the triangulation of the solid torus (for the time being, $x_{i}$ is just a formal variable);
\item two opposite edges carrying the same dihedral angle $y_{i}$ and identified to just one edge, of slope $\eta$, in the triangulation (similarly, $y_{i}$ is formal);
\item two opposite edges which carry the same (formal) dihedral angle $z_{i}$, and which coincide with the edges of slope $\zeta$ and $\zeta'$ in the triangulation.
\end{itemize}
As in any angle structure, the relationship $x_{i}+y_{i}+z_{i}=\pi$ must hold between the formal variables.

The vertices of the hexagon $H_{i-1}$ (resp. $H_i$) are the links of edges of slopes $\xi,\eta,\zeta$ (resp. $\xi,\eta,\zeta'$). We can write these labels $\xi,\eta,\zeta,\zeta'$ at the vertices of $H_{i-1}$ and $H_i$: see Figure \ref{fig:hexa} (left).

\begin{observation} \label{obs:angles}
By construction, the vertex of the hexagon $H_{i-1}$ labelled $\zeta$ has an interior angle of $z_{i}$, while the vertex of hexagon $H_i$ labelled $\zeta'$ has an interior angle of $2\pi-z_{i}$. This comes from the fact that the boundary of the tetrahedron $\Delta_{i}$ is exactly the union of the two pleated punctured tori $\tau_{i-1}$ and $\tau_i$ (with vertex links $H_{i-1}, H_i$).
\end{observation}

As a consequence, we can determine the three angles of the hexagon $H_i$ (each angle occurs, by central symmetry of $H_i$): 
\begin{equation} \label{eq:3angles} 2\pi-z_{i}\hspace{8pt} ; \hspace{8pt} z_{i+1} \hspace{8pt} ; \hspace{8pt} z_{i}-z_{i+1}\:.\end{equation}
Indeed, the first two of these numbers are given by Observation \ref{obs:angles} (shifting indices by one for $z_{i+1}$); the third is given by the property that the six angles of $H_i$ should add up to $4\pi$. See Figure \ref{fig:hexa}, (right).

\begin{figure}[h!] \centering
\psfrag{zp}{$z_{2}$}
\psfrag{zm}{$z_{1}$}
\psfrag{H0}{$H_0$}
\psfrag{H1}{$H_1$}
\psfrag{H2}{$H_2$}
\psfrag{x}{$\xi$}
\psfrag{y}{$\eta$}
\psfrag{z}{$\zeta$}
\psfrag{zz}{$\zeta'$}
\includegraphics{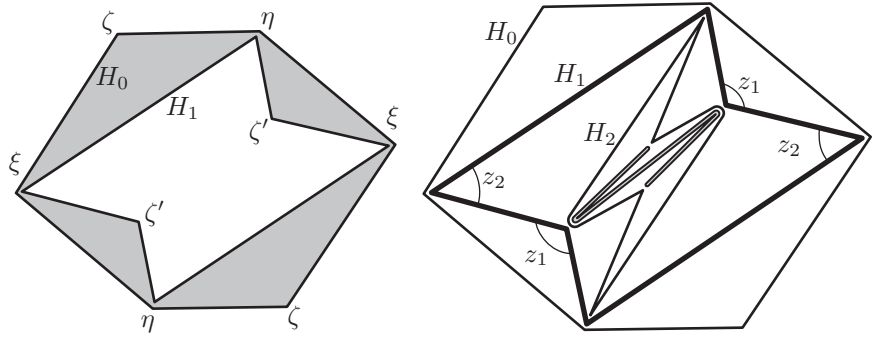}
\caption{Left: two consecutive hexagons $H_0, H_1$ in the cusp link, with vertices labelled by elements of $\mathbb{P}^1\mathbb{Q}$. The four (similar) grey triangles are the vertex links of the ideal tetrahedron $\Delta_{1}$. Right: the full sequence of hexagons $H_0,\dots,H_3$, where $H_3$ is collapsed to a broken line of $3$ segments. The angles $z_{1}$ and $z_{2}$ of the tetrahedra $\Delta_{1}$ and $\Delta_{2}$ are marked; together they determine the interior angles (\ref{eq:3angles}) of $H_1$.} \label{fig:hexa} 
\end{figure}

We can in turn write the numbers (\ref{eq:3angles}) in the corners of the Farey triangle $T_i$: namely, $2\pi-z_{i}$ is in the corner opposite the Farey edge $T_{i-1}\cap T_i$; similarly $z_{i+1}$ is in the corner opposite the Farey edge $T_i\cap T_{i+1}$; and $z_{i}-z_{i+1}$ is in the third corner, at the Farey vertex $T_{i-1}\cap T_{i+1}$. See Figure \ref{fig:rlrl}.

The above operation can be performed for all indices $i\in \{1,\dots,N-2\}$. For $i=N-1$, there is no tetrahedron ``$\Delta_{N}$''; hence, a priori, no parameter $z_{N}$. However, if $m'\in \mathbb{P}^1\mathbb{Q}$ is the vertex of the Farey triangle $T_{N-1}$ opposite the Farey edge $T_{N-1}\cap T_N$ in $T_{N-1}$, then the interior angle of the (collapsed) hexagon $H_{N-1}$ at the vertex labelled $m'$ is precisely $0$, by definition of our folding of the pleated surface $\tau_{N-1}$ onto itself. This folding thus corresponds to asking that $$z_{N}=0~.$$ Under this convention, the other angles of the collapsed hexagon $H_{N-1}$ are then given by the same formulas (\ref{eq:3angles}), with $i=N-1$.

Finally, we perform an analogous construction at $i=0$ (it follows from our assumptions that $H_0$ is convex, with angles $\pi-\theta_p, \pi-\theta_q, \pi-\theta_r$). There is no tetrahedron ``$\Delta_{0}$''; hence, a priori, no parameter $z_{0}$.
%if $m'\in \mathbb{P}^1\mathbb{Q}$ is the vertex of the Farey triangle $T_{N-1}$ opposite the Farey edge $T_{N-1}\cap T_N$ in $T_{N-1}$, then 
However, the interior angle of $H_0$ at the vertex labelled $r$ is $\pi-\theta_r$, which entails $z_{1}=\pi-\theta_r$. Similarly, the interior angle of $H_0$ at the vertex labelled $q$ is $\pi-\theta_q$, which entails $z_{0}=2\pi-(\pi-\theta_q)=\pi+\theta_q$. To summarize,

\begin{proposition} Under the full set of assumptions
\begin{equation} \label{eq:naildown} 
\begin{array}{rrrrrrr}
(&z_{0}~,& z_{1}~,& z_{2}~,&\dots~,& z_{N-1}~,&z_{N}~) \\
=~(&\pi+\theta_{q}~,& \pi-\theta_r~,& z_{2}~,&\dots~,& z_{N-1}~,&0~) \end{array} \end{equation} (where the values of $(z_{2},\dots, z_{N-1})$ remain to be chosen), the angles of the hexagons $\{H_i\}_{0\leq i \leq N-1}$ given by (\ref{eq:3angles}) define all the $(\theta_p, \theta_q, \theta_r)$--angle structures. \qed
\end{proposition}

To get angle structures, we must only choose the $z_{2},\dots, z_{N-1}$ in the interval $(0,\pi)$ so that all dihedral angles of $\Delta_{i}$ are positive for $1\leq i \leq N-1$, which we do now.

Denote by $\xi$ (resp. $\eta$) the right (resp. left) end of the Farey edge $T_{i-1}\cap T_i$. By construction, $x_{i}$ (resp. $y_{i}$) is half the difference between the angles of hexagons $H_{i-1}$ and $H_i$ at the vertex labelled $\xi$ (resp. $\eta$) in the cusp link, i.e. half the difference between the numbers written in the $\xi$--corner (resp. the $\eta$--corner) of the Farey triangles $T_{i-1}$ and $T_i$ in the Farey diagram. (The factor one--half comes from the identification of pairs of opposite edges in the ideal tetrahedron $\Delta_{i}$.) In Figure \ref{fig:rlrl} we show what these numbers are, according to whether the line $\ell$ makes Rights or Lefts at the Farey triangles $T_{i-1}$ and $T_i$: we use only (\ref{eq:3angles}) and the shorthand 
\begin{equation} \label{eq:abczi}
(a,b,c):=(z_{i-1}~,~z_{i}~,~z_{i+1})~. \end{equation}

\begin{figure}[h!] \centering
\psfrag{x}{$\xi$}
\psfrag{y}{$\eta$}
\psfrag{r}{$R$}
\psfrag{l}{$L$}
\psfrag{b}{$b$}
\psfrag{pb}{$2\pi\text{-}b$}
\psfrag{pa}{$2\pi\text{-}a$}
\psfrag{c}{$c$}
\psfrag{ab}{$a\text{-}b$}
\psfrag{bc}{$b\text{-}c$}
\includegraphics[scale=.7]{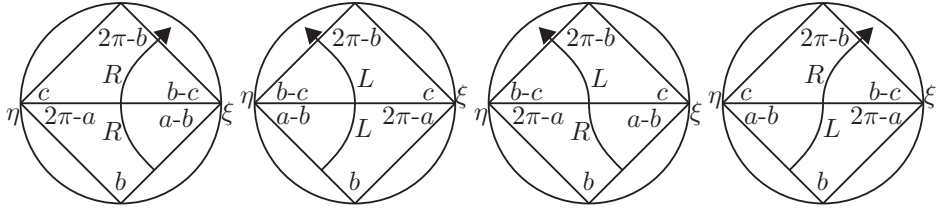}
\caption{The Farey triangles $T_{i-1}$ (lower) and $T_i$ (upper), with corner labels.} \label{fig:rlrl} 
\end{figure}

It follows that the values of $x_{i}$ and $y_{i}$ in terms of the $z_{i}$ are given by Table (\ref{tab:xiyi}) --- in the first line of the table, we recall the nature of the tetrahedron (or cell) $\Delta_{i}$, and the natural positions of $a,b,c$, interspersed with the letters of the word $\Omega$.

\setstretch{2.6} 
\begin{equation}
\label{tab:xiyi}
\begin{array}{c|cccc}
\underset{\text{Cell $\Delta_{i}$ is...}}{\overset{z_{i\text{-}\!1}\hspace{6pt}z_{i}\hspace{6pt}z_{i+\!1}}{\underbrace{\Omega_{i-1} \: , \: \Omega_i}}}
& \underset{\text{Non--hinge}}{\overset{a\hspace{12pt}b\hspace{12pt}c}{\underbrace{R\hspace{3pt},\hspace{3pt}R}}} 
& \underset{\text{Non--hinge}}{\overset{a\hspace{12pt}b\hspace{12pt}c}{\underbrace{L\hspace{3pt},\hspace{3pt}L}}} 
& \underset{\text{Hinge}}{\overset{a\hspace{12pt}b\hspace{12pt}c}
{\underbrace{R\hspace{3pt},\hspace{3pt}L}}} 
& \underset{\text{Hinge}}{\overset{a\hspace{12pt}b\hspace{12pt}c}
{\underbrace{L\hspace{3pt},\hspace{3pt}R}}} 
\\ \hline
x_{i} &\displaystyle{\frac{a-2b+c}{2}} & \pi-\displaystyle{\frac{a+c}{2}}& \displaystyle{\frac{a-b-c}{2}} & \pi-\displaystyle{\frac{a+b-c}{2}} \\ 
y_{i} & \pi-\displaystyle{\frac{a+c}{2}}& \displaystyle{\frac{a-2b+c}{2}}& \pi-\displaystyle{\frac{a+b-c}{2}} & \displaystyle{\frac{a-b-c}{2}} \\  
z_{i} & b & b & b & b \end{array}
\end{equation}
\setstretch{1.0}

From Table (\ref{tab:xiyi}), we can read off the condition for all $x_{i}$ and $y_{i}$ and $z_{i}$ to be positive. Still using the notation $(a,b,c)=(z_{i-1}~,~z_{i}~,~z_{i+1})$, these conditions are
\begin{equation} \label{racohi} \left \{ \begin{array}{ll}
\bullet \hspace{8pt} 
%z_{i-1}>z_{i}+z_{i+1} 
a>b+c
& \text{if $\Delta_{i}$ is a hinge cell (\emph{hinge condition});} \\
\bullet \hspace{8pt} 
%z_{i-1}+z_{i+1}>2z_{i} 
a+c>2b
& \text{if $\Delta_{i}$ is not a hinge (\emph{convexity condition});} \\
\bullet \hspace{8pt} 
0<z_{i}<\pi 
%0<b<\pi
& \text{for all $2\leq i \leq N-1$ (\emph{range condition});} \\
\bullet \hspace{8pt} z_{2}<\pi-\theta_q & \text{(follows from the case $i=1$, a non--hinge index).}
\end{array} \right . \end{equation}
The last condition is needed for $\pi-\frac{z_{0}+z_{2}}{2}$ (namely, $x_{1}$ or $y_{1}$) to be positive, because $z_{0}=\pi+\theta_q$ (unlike other $z_{i}$) is larger than $\pi$. Note that by (\ref{eq:naildown}), the convexity condition at $i=1$ also implies $z_{2}>\pi-\theta_q-2\theta_r$. This is compatible with the last condition of (\ref{racohi}) since $\theta_r>0$ by (\ref{eq:thetas}).

\begin{itemize}
\item {\bf Case 1: none of the $\Delta_{i}$ are hinge cells.}
In this case, we are reduced to finding a sequence of the form (\ref{eq:naildown}) that is convex, decreasing, and satisfies $z_{2}<\pi-\theta_q$. This is clearly possible if and only if
$$\begin{array}{rrcl}
&(\pi+\theta_q)-N(\theta_r+\theta_q)&<&0~, \\
\text{i.e.}& (N-1)\theta_q + N\theta_r &>&\pi~, \\
\text{i.e.}& \theta_p + N \theta_q + (N+1)\theta_r &>& 2\pi~, \end{array}$$
where the last line follows from (\ref{eq:thetas}). It is easy to check that under the normalization $(p,q)=(\infty,0)$ and $r\in \{+1,-1\}$ (one of which can be assumed up to applying an element of $PSL_2(\mathbb{Z})$), the slope $m\in\mathbb{P}^1\mathbb{Q}$ is, up to sign, the integer $N$: indeed, all the letters of the word $\Omega$ are equal and the Farey triangle $T_i$ has vertices $\infty,i,i-1$ if $r=1$ (and $\infty,-i,-i+1$ if $r=-1$). The last line of the computation above thus becomes $$ (m\wedge p)\theta_p + (m\wedge q)\theta_q + (m\wedge r)\theta_r > 2\pi~,$$ proving Proposition \ref{prop:pqr} in this case.

\item {\bf Case 2: some $\Delta_{i}$ are hinge cells.} By Remark \ref{rem:zebre}, the inequality of Proposition \ref{prop:pqr} is vacuous in this case. Let us therefore just construct a sequence of the form (\ref{eq:naildown}) that satisfies (\ref{racohi}). Let $h \in \{2,3,\dots,N-1\}$ be the smallest hinge index. We can easily choose a strictly convex, positive, decreasing sequence 
$$ \begin{array}{rrrrrrr}
(&z_{0}~,& z_{1}~,& z_{2}~,&\dots~,& z_{h-1}~,&z_{h}~) \\
=~(&\pi+\theta_{q}~,& \pi-\theta_r~,& z_{2}~,&\dots~,& z_{h-1}~,&z_{h}~) \end{array} $$
satisfying $z_{2}<\pi-\theta_q$. We construct the rest of the sequence $(z_{i})$ backwards, descending from the index $i=N$ down to $i=h$. First set $z'_{N}=0$ and $z'_{N-1}=1$. For each $i$ such that $N-2\geq i \geq h+1$, pick (inductively) a value of $z'_{i}$ such that $(a,b,c):=(z'_{i},z'_{i+1},z'_{i+2})$ satisfies the concavity or hinge condition of (\ref{racohi}), according to whether $\Delta_{i+1}$ is a hinge cell or not (for example, $z'_{i}=3z'_{i+1}$ will do). The sequence $(z'_{h+1},\dots,z'_{N-1})$ is clearly positive and decreasing. We then set $$z_{i}:=\varepsilon z'_{i} ~ \text{ for all }~ h+1\leq i \leq N~:$$ it is immediate to check that the hinge condition ``$a>b+c$'' of (\ref{racohi}) is verified by the triple $(a,b,c)=(z_{h-1},z_{h},z_{h+1})$ as soon as $$0<\varepsilon< \frac{z_{h-1}-z_{h}}{z'_{h+1}}.$$ Thus, by choosing such an $\varepsilon$, we have found a sequence $(z_{i})$ of the form (\ref{eq:naildown}).\end{itemize} 
Proposition \ref{prop:pqr} is proved.
\end{proof}

\subsection{Volume maximization}
\label{sec:maximize}
Denote by (\ref{racohi}') the system (\ref{racohi}) in which all strong inequalities have been replaced by weak ones, and let $W$ denote the compact polyhedron of solutions $(z_{i})$ of the form (\ref{eq:naildown}) to the system (\ref{racohi}') --- so the interior of $W$ is the \emph{space of angle structures}. The volume functional $\mathcal{V}:W\rightarrow \mathbb{R}^+$ associates to every point $z$ of $W$ the sum of the volumes of the ideal tetrahedra $\Delta_{i}$ with non-negative angles $x_{i}, y_{i},z_{i}$ given by Table \ref{tab:xiyi}. 

Suppose that $\theta_p, \theta_q, \theta_r$ satisfy (\ref{eq:thetas}) and the inequality of Proposition \ref{prop:pqr} (hence $W\neq \emptyset$). We henceforth assume that the point $z=(z_{i})\in W$ realizes the maximum of $\mathcal{V}$ over $W$, and we aim to prove 

\begin{proposition} \label{winterior}
The point $z$ is a solution of (\ref{racohi}), not just (\ref{racohi}') --- i.e., all $\Delta_{i}$ have only positive angles.
\end{proposition}
\begin{proof}
Observe that the sequence $(z_{0},\dots,z_{N})$ is non-negative and non-increasing. This follows from (\ref{racohi}') by an immediate downward induction (starting at $z_{N}$).

By Proposition \ref{prop:degeneracies}, we know that if $\Delta_{i}$ is \emph{flat}, i.e. has a vanishing dihedral angle, then its triple of angles is of the form $(0,0,\pi)$, up to permutation. Thus, by Table (\ref{tab:xiyi}), $\Delta_{i}$ is flat exactly when $z_{i}\in\{0,\pi\}$. By monotonicity, since $z_{1}=\pi-\theta_r<\pi$, the only flat tetrahedra $\Delta_{i}$ actually satisfy $z_{i}=0$. Still by monotonicity, it then follows that $z_{i+1}=0$ as well. Let $i$ be the \emph{smallest} index such that $z_{i}=0$. An easy discussion, using Table (\ref{tab:xiyi}), shows that the only possible value of $z_{i-1}$ that implies $\{x_{i}\, y_{i}\}=\{0,\pi\}$ is $z_{i-1}=2\pi$ (recall here the $a$-$b$-$c$--notation \ref{eq:abczi}). This is impossible: only $z_{0}=\pi+\theta_q$ is allowed to be larger than $\pi$, but we have $\theta_q<\pi$ by (\ref{eq:thetas}).
\end{proof}

\begin{corollary} \label{cor:completeness}
The point $z$ defines a complete hyperbolic structure on the punctured solid torus $X=\Delta_{1}\cup\dots\cup \Delta_{N-1}$, with exterior dihedral angles $\theta_p, \theta_q, \theta_r$ on $\partial X$.
\end{corollary}
\begin{proof}
By Theorem \ref{thm:rivin}, this follows from the fact that $z$ is critical for the volume functional $\mathcal{V}:W\rightarrow \mathbb{R}$.
%An alternative proof would closely follow that of \cite[Lemma 6.2]{gufu}: to each interior edge $E$ of $X$ is associated a certain line $L_E$ in the tangent space $T_z W$, such that the vanishing of the derivative of $\mathcal{V}$ along $L_E$ expresses the fact that the hyperbolic metric near $E$ is complete.
\end{proof}

\section{Handedness}
\label{sec:handedness}
%%%%%%%%%%%%%%%%%%%%%%%%
%%                    %%
%%   handedness.tex   %%
%%                    %%
%%%%%%%%%%%%%%%%%%%%%%%%

In this section, we discuss the \emph{handednesses} of certain elements in the fundamental group of the (complete, hyperbolic) punctured solid torus $X$. These results will be useful in establishing the inequalities leading to Theorem \ref{thm:main} (which is proved in the next section).

\begin{definition} \label{def:handedness}
For any $g\in GL_2(\mathbb{C})$, define the handedness of $g$ by $${\sf hand}\,(g):=\frac{(\text{Tr}\, g)^2}{\text{Det}\, g}~.$$ Note that ${\sf hand}\,(g)={\sf hand}\,(g^{-1})={\sf hand}\,(rg)$ for all $r\in \mathbb{C}^*$. Therefore, ${\sf hand}$ factors through a map $PSL_2(\mathbb{C})\rightarrow \mathbb{C}$, also noted ${\sf hand}$. Call a loxodromy of $\mathbb{H}^3$ \emph{left--handed} (resp. \emph{right--handed}) when it is conjugate to $z\mapsto \alpha z$ with $|\alpha|>1$ and $\text{Im}\,(\alpha)>0$ (resp. $|\alpha|>1$ and $\text{Im}\,(\alpha)<0$). Left--handed loxodromies are ``corkscrew'' motions, the motion of a dancer who jumps upwards while spinning to his left. It is easy to check that the M\"obius transformation associated to $g$ is left-- (resp. right--) handed if and only if $\text{Im}\,({\sf hand}\,(g))$ is positive (resp. negative).
\end{definition}

Let $U$ be a universal cover of the solid torus $X=\bigcup_{i=1}^{N-1} \Delta_{i}$. Since $U$ is a complete hyperbolic manifold with locally convex boundary, the developing map $U\rightarrow \mathbb{H}^3$ is an embedding. Thus $U\subset \mathbb{H}^3$ is the convex hull in $\mathbb{H}^3$ of the orbit of an ideal point $v$ under a certain loxodromy $$\varphi \in \text{Isom}^+(\mathbb{H}^3)\simeq PSL_2(\mathbb{C})$$ (typically extremely short, corresponding to the core curve of the solid torus). Make the attractive (resp. repulsive) fixed point of $\varphi$ coincide with the North pole $P^+$ (resp. the South pole $P^-$) of $\mathbb{S}^2\simeq\partial_{\infty}\mathbb{H}^3$; assume that $v$ lies on the Equator at longitude $0$, and orient the Equator along increasing longitudes. As a cover of the space $X$ which is triangulated, $U$ comes with a natural, $\varphi$--invariant decomposition into ideal tetrahedra. 

The projection with respect to the center of Poincar\'e's ball model sends $\partial U$ homeomorphically to $\mathbb{S}^2\smallsetminus\{P^+,P^-\}\smallsetminus\{\varphi^n(v)\}_{n\in\mathbb{Z}}$. For each edge $vv'$ of $\partial U$ (between ideal points $v,v' \in \mathbb{S}^2$), this projection sends $vv'$ to the short great--circle arc $\overset{\frown}{vv'}$ in $\mathbb{S}^2$. If $vv''$ is another edge of $\partial U$, this allows us to speak about the \emph{angle} $\widehat{v'vv''}\in (-\pi,\pi]$ between $v'$ and $v''$, as seen from $v$ (i.e. in $T_v\mathbb{S}^2$).

The punctured torus $\tau_0=\partial U /\varphi$ has three ideal edges, each endowed with a positive dihedral angle. Therefore the ideal vertex $v$ of $U$ is connected to six other vertices of $U$ by edges of $\partial U$, and there is a natural cyclic order on these six vertices. The equatorial plane intersects $\partial U$ along a broken line $J$ from $v$ to $v$ which is properly embedded in $\partial U$ (with ideal endpoints). We can orient $J$ along increasing longitudes. 

\begin{figure}[h!] \centering
\psfrag{p}{$v$}
\psfrag{p1}{$v_1$}
\psfrag{p2}{$v_2$}
\psfrag{p3}{$v_3$}
\psfrag{p4}{$v_4$}
\psfrag{p5}{$v_5$}
\psfrag{p6}{$v_6$}
\psfrag{q}{Equator}
\psfrag{N}{$P^+$}
\psfrag{S}{$P^-$}
\psfrag{f}{$\varphi(v)$}
\psfrag{H}{$\mathbb{S}^2\simeq \partial_{\infty} \mathbb{H}^3$}
\includegraphics[scale=.7]{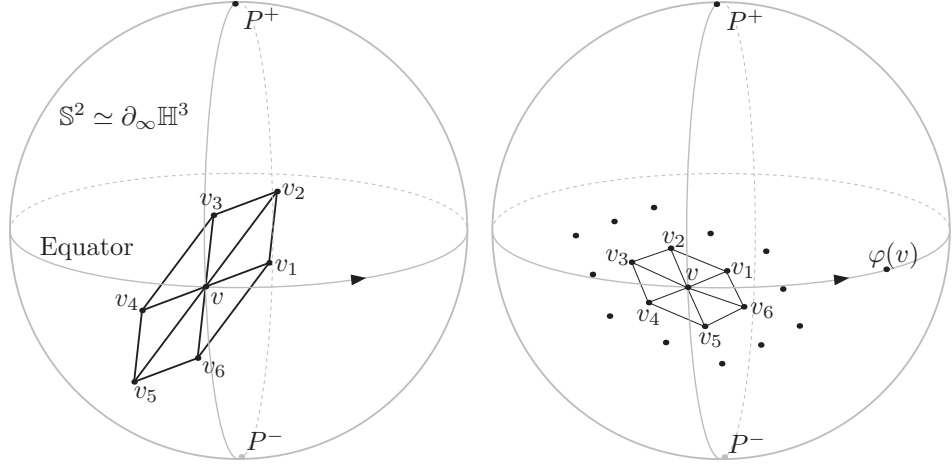}
\caption{Left: one cannot have $l_6 \leq 0 < l_1$. Right: the actual situation (only some ideal vertices of $U$ are shown).}
\label{fig:globe}
\end{figure}

\begin{definition} \label{def:sixneighbors}
Let $v_1,\dots v_6$ (with indices seen modulo $6$) denote the six neighbors of $v$ that are met, in that order, when turning counterclockwise around $v$, starting in the direction of the initial segment of $J$. For each $i$ in $\mathbb{Z}/6\mathbb{Z}$, there is an integer $n_i\in \mathbb{Z}$ such that $\varphi^{n_i}$ sends the following points to one another: $$\begin{array}{cc} & ~ v_{i+2} \mapsto v_{i+1} \\ \varphi^{n_i}~: & v_{i\pm 3} \longmapsto v \longmapsto v_i ~ \\ & ~ v_{i-2} \mapsto v_{i-1}. \end{array} $$ Of course, $n_i=-n_{i+3}$. See Figure \ref{fig:globe}.
\end{definition}

\begin{claim} \label{claim:longitudes}
The longitudes $l_1, l_6$ of $v_1$ and $v_6$ are both in $(0,\pi)$. The latitude of $v_1$ (resp. $v_6$) is positive (resp. negative).
\end{claim}

\begin{proof}
Since a half-turn around $v$ sends each $v_i$ to $v_{i+3}$, no angle $\widehat{v_{i-1} v v_i}$ in the tangent space to $\mathbb{S}^2$ at $v$ can exceed (or even reach) the value $\pi$: taking $i=1$, this proves the statement about latitudes. Therefore $v_1$ (resp. $v_6$) lies above (resp. below) the equatorial plane, and it also follows that $n_6<0<n_1$. 

Let $l_i \in (-\pi,\pi]$ denote the longitude of $v_i$: clearly, $l_i<\pi$ since no edge of $\partial U$ can cross the North-South axis. The longitudes $l_1$ and $l_6$ cannot be both nonpositive, otherwise $\widehat{v_6 v v_1}\geq \pi$. Therefore, assume $$l_6 \leq 0 < l_1$$ and aim at a contradiction. 

Note that on $\mathbb{S}^2$, for each $n>0$, the transformation $\varphi^n$ increases latitudes, and adds a constant angle to all longitudes (modulo $2\pi$). Recall the relationships $v_3=\varphi^{-n_6}(v)$ and $v_2=\varphi^{n_1}(v_3)=\varphi^{-n_6}(v_1)$: they imply that $v_2$ has highest \emph{latitude} among $v_1,v_2,v_3$ (all three latitudes being positive). They also imply $l_2\equiv l_1-l_6~[\text{mod } 2\pi]$: but $l_2$ cannot belong to $(\pi,l_1+\pi)+2\pi\mathbb{Z}$ since the ideal triangle $vv_1v_2\subset \partial U$ cannot meet the North-South axis. Therefore, $l_2=l_1-l_6=l_1+l_3$ belongs to $(l_1,\pi)$, and the point $v_2$ also has the largest \emph{longitude} among $v_1,v_2,v_3$, possibly tying with $v_1$ (and all three longitudes belong to $[0,\pi)$).

It follows that the hyperbolic line $vv_2$ comes closer to the North-South axis than the hyperbolic line $v_1v_3$, which contradicts the convexity of $U$ near the edge $vv_2$: absurd. See Figure \ref{fig:globe}.
\end{proof}

\begin{remark}
Claim \ref{claim:longitudes} implies that $\varphi^{\pm n_1}$ and $\varphi^{\pm n_6}$ are, respectively, left-- and right-handed.
\end{remark}

\medskip

%Orient the edges $vv1$ and $vv_6$ of the infinite polyhedron $U$ along increasing longitudes, or, equivalently, away from the vertex $v$: we write $\overrightarrow{vv_1},\overrightarrow{vv_6}$ for these oriented edges, and denote by $t$ the ideal triangle $vv_1v_6 \subset \partial U$ that connects them.

Recall the sequence of Farey triangles $pqr=T_0,T_1,\dots,T_N=mst$. All $T_i$ for $i\geq 1$ have their vertices in the arc $\overset{\frown}{pq} \subset \mathbb{P}^1\mathbb{Q}$ that does not contain $r$ (in particular, the meridinal slope $m$ belongs to that arc). For every $i\in \{1,2,\dots,N\}$ and every vertex $x$ of the Farey triangle $T_i$, we can draw a properly embedded intrinsic geodesic $g_{x}$ of slope $x$ in the punctured torus $\partial U/\varphi$: this $g_{x}$ has a lift $\widehat{g_{x}}\subset \partial U$ that connects the ideal vertex $v$ 
to some $\varphi$--iterate of $v$, and whose initial (ideal) segment is contained in the ideal triangle $vv_1v_6$ of $\partial U$. We orient $\widehat{g_{x}}$ from $v$ to its other end. 

\begin{definition}
When $x\in\mathbb{P}^1\mathbb{Q}$ is a vertex of some Farey triangle $T_i$ as above, define $\nu_x \in \mathbb{Z}$ as the integer such that the oriented curve $\widehat{g_x}$ runs from the ideal vertex $v$ to $\varphi^{\nu_x}(v)$. 

We also define $\lambda_{x} \in \mathbb{R}$ as the integral of the longitude $1$--form in $\mathbb{S}^2 \smallsetminus \{P^+,P^-\}$ along the closure of $\pi(\widehat{g_{x}})$, where $\pi:\partial U \rightarrow \mathbb{S}^2$ is the central projection.
\end{definition}

\begin{proposition}\label{prop:fibonacci}
Suppose $1\leq i \leq N-1$ so that  $T_i=abc$ and $T_{i+1}=bcd$ are two consecutive Farey triangles. Then $\nu_d=\nu_b+\nu_c$ and $\lambda_d=\lambda_b+\lambda_c$. 

Moreover, if $x \in \mathbb{P}^1\mathbb{Q}$ is a vertex of $T_i$ for some $i\in\{1,\dots,N\}$, then $0<\lambda_{x}\leq 2\pi$, with equality (for the upper bound) if and only if $x$ is the meridinal slope $m$.
\end{proposition}
\begin{proof}
Consider the ideal quadrilateral $Q:=(\partial U/\varphi) \smallsetminus (g_b \cup g_c)$. The orientations of $g_b$ and $g_c$ induce orientations on the four edges of (the metric completion of) $Q$. Observe that $g_d$ runs diagonally across $Q$, from the vertex with two outgoing edges, to the vertex with two incoming edges: as a result, the closure of $\pi(\widehat{g_d})$ in $\mathbb{S}^2 \smallsetminus \{P^+,P^-\}$ is isotopic with endpoints fixed to the closure of $$\pi \left ( \widehat{g_b} \cup \varphi^{\nu_b}(\widehat{g_c}) \right ) ~\text{ or, indifferently, of }~\pi \left ( \widehat{g_c} \cup \varphi^{\nu_c}(\widehat{g_b}) \right ).$$ The exponent identity $\nu_d=\nu_b+\nu_c$ follows and, since $\varphi$ increases longitudes by a constant, so does the longitude identity $\lambda_d=\lambda_b+\lambda_c$. 

By Claim \ref{claim:longitudes}, we have $\lambda_p, \lambda_q \in (0,\pi)$, so an immediate upward induction on $i$ now implies $\lambda_{x}>0$ for each vertex $x$ of $T_i$ (with $1\leq i \leq N$). But $\lambda_m=\pm 2\pi$, because the meridian curve $\widehat{g_m}$ runs exactly once around the infinite polyhedron $U$: therefore, $\lambda_m=2\pi$. Downward induction on $i$ finally yields $\lambda_{x}<2\pi$ for $x\neq m$.
\end{proof}

\begin{proposition} \label{prop:handednesses}
Suppose $1\leq i \leq N-1$. Let $x \in \mathbb{P}^1\mathbb{Q}$ be the Farey vertex common to $T_{i-1},T_i,T_{i+1}$. Then,
\begin{enumerate}[(i)]
\item one has $\lambda_{x} \in (0,\pi);$
\item if the Farey triangle $T_i$ carries an $L$ (resp. an $R$), then $\nu_x>0$ (resp. $\nu_x<0$);
\item if $T_i$ carries an $L$ (resp. an $R$), then $\varphi^{\nu_x}$ is left--handed (resp. right--handed).
\end{enumerate}
\end{proposition}

\begin{proof}
We name the vertices of the Farey triangles so that $T_i=xyz$ and $T_{i+1}=xzt$. By Proposition \ref{prop:fibonacci}, one has $\lambda_z=\lambda_x+\lambda_y$ and $2\pi\geq \lambda_t=\lambda_x+\lambda_z=2\lambda_x+\lambda_y$. Since $\lambda_x, \lambda_y>0$, this yields (i).

Assertion (ii) follows from the following claim: if $l_i,r_i \in \mathbb{P}^1\mathbb{Q}$ are the left and right endpoints of the Farey edge $T_{i-1}\cap T_i$ (for the transverse orientation towards $m$), then $\nu_{r_i}<0<\nu_{l_i}$. This is clearly true for $i=1$ (in that case, $\nu_{l_i}=n_1$ and $\nu_{r_i}=n_6$, in the notation of Definition \ref{def:sixneighbors}). For $i>1$, observe that 
\begin{itemize}
\item one has $\nu_m=0$ because the curve $\widehat{g_m}$ is a closed curve around the ideal polyhedron $U$;
\item by Proposition \ref{prop:fibonacci}, the number $\nu_m$ is always a linear combination of $\nu_{l_i}$ and $\nu_{r_i}$ with positive integer coefficients;
\item one has $\nu_{l_i}\neq 0$ and $\nu_{r_i}\neq 0$ because the curves $\widehat{g_{l_i}}$ and $\widehat{g_{r_i}}$ are \emph{not} closed curves in $\partial U$.
\end{itemize}
These observations put together imply $\nu_{r_i}<0<\nu_{l_i}$ or $\nu_{l_i}<0<\nu_{r_i}$. The first is clearly the case by induction on $i$, because one always has $l_i=l_{i+1}$ (resp. $r_i=r_{i+1}$) if the Farey triangle $T_i$ carries an $L$ (resp. an $R$).

Assertion (iii) is an immediate consequence of (i)--(ii). 
\end{proof}

%%%%%%%%%%%%%%%%%%%%%%%%%%%%%%

\section{Canonical decomposition of a generic Dehn filling} 
\label{sec:voronoi}
%%%%%%%%%%%%%%%%%%%%%
%%                 %%
%%   voronoi.tex   %%
%%                 %%
%%%%%%%%%%%%%%%%%%%%%

In this section we prove Theorem \ref{thm:main}: to show that a given triangulation is Delaunay (or geometrically canonical), we essentially must prove a certain number of inequalities, which will boil down to statements of handedness as given by Proposition \ref{prop:handednesses}. 

Consider a hyperbolic manifold $M$ with $k\geq 2$ cusps, endowed with horoball neighborhoods, such that the genericity assumption of Theorem \ref{thm:main} is satisfied. Let $\mathcal{D}$ denote the canonical triangulation of $M$. We assume that $H_k$, the horoball neighborhood of the $k$--th cusp $c_k$, has much smaller volume than all other $H_i$.

\subsection{A generic small cusp}
First we prove that $\mathcal{D}$ contains exactly two ideal tetrahedra $\Delta,\Delta'$ that have a vertex in $c_k$.

Consider a universal covering $\pi:\mathbb{H}^3\rightarrow M$ such that (in the upper half--space model) the point at infinity lies above the cusp $c_k$. Let $\Lambda$ be the rank--$2$ lattice of deck transformations of the form $z\mapsto z+\lambda$. Let $\{\eta_i\}_{i\in I}$ be the collection of all horoballs of $\mathbb{H}^3$ lying above some $H_i$ with $i<k$ (the $\eta_i$ are Euclidean balls tangent to the boundary $\mathbb{C}$ of the model half--space.) By the genericity assumption of Theorem \ref{thm:main}, there is a unique shortest path in $M$ from $H_k$ to $\bigcup_{i=1}^{k-1}H_i$: therefore the largest $\eta_i$ (for the Euclidean metric) is unique modulo $\Lambda$.

We can assimilate $\Lambda$ to a lattice of $\mathbb{C}$, and assume that the largest $\eta_i$'s are centered exactly at the points of $\Lambda$.

The Delaunay decomposition $D_{\Lambda}$ of $\mathbb{C}$ with respect to the vertex set $\Lambda$ consists either of isometric rectangles (all belonging to the same $\Lambda$--orbit), or of isometric triangles (belonging to two $\Lambda$--orbits) with strictly acute angles. We claim that the latter is the case: indeed, let $P\subset \mathbb{C}$ be a convex polygon of $D_{\Lambda}$: the vertices of $P$, which are points of $\Lambda$, are on the boundary of a disk that contains no other points of $\Lambda$. Using the fact that the horoball $\eta_{\infty}$ centered at infinity stays very high above $\mathbb{C}$ in the half--space model (because $H_k$ has very small volume), it is easy to construct a ball of $\mathbb{H}^3$ that is tangent to the horoballs $\eta_i$ centered at the vertices of $P$, disjoint from all other $\eta_i$, and tangent to the horoball $\eta_{\infty}$. The center of this ball is a vertex of the Ford domain. Hence, there exists a cell of the Delaunay decomposition $\mathcal{D}$ of $M$ (more precisely, a lift $\widehat{\Delta}$ of such a cell to $\mathbb{H}^3$) whose vertices are precisely $\infty$ and the vertices of $P$. By the genericity assumption, $\widehat{\Delta}$ must be an ideal tetrahedron, so $P$ is a triangle, and has strictly acute angles. The two (isometric) $\Lambda$--orbits of triangles in the Delaunay decomposition $D_{\Lambda}$ of $\mathbb{C}$ correspond to two ideal tetrahedra $\Delta,\Delta'$ in $\mathcal{D}$. Note that $\Delta\cup \Delta'$ is a neighborhood of the cusp $c_k$.

The space $T=\partial (\Delta\cup\Delta')\subset M$ is the quotient by $\Lambda$ of the union of all ideal triangles of $\mathbb{H}^3$ that project vertically to triangles of $D_{\Lambda}$ (contained in $\mathbb{C}$): therefore, $T$ is a hyperbolic once--punctured torus bent along three lines, and its interior dihedral angles are twice those of $\Delta$ (or $\Delta'$).

\subsection{Triangulation of the Dehn filling}

It is well--known that almost all (hyperbolic) Dehn fillings $M_s$ of $M$ at the cusp $c_k$ admit a \emph{spun} decomposition $\mathcal{D}_s^{\text{spun}}$ into ideal, positively--oriented tetrahedra: namely, $\mathcal{D}_s^{\text{spun}}$ is obtained from $\mathcal{D}$ by letting the tips of $\Delta$ and $\Delta'$ (formerly in $c_k$) spin asymptotically along the geodesic core of the filling solid torus of $M_s$ --- actually, there are two such spun triangulations, spinning in opposite directions (see e.g. \cite{thurstonnotes}, Chap. V). Moreover, the cross--ratios of the tetrahedra of $\mathcal{D}_s^{\text{spun}}$ become (uniformly) close to those of $\mathcal{D}$ as the slope $s$ goes to infinity, i.e. gets more and more complicated. In particular, the punctured torus $T$, equal to the union of the bases of $\Delta$ and $\Delta'$, is still embedded in $M_s$, with bending angles close to those in $M$.

Therefore, we can remove the solid torus $\Delta\cup\Delta'$ from the spun triangulation of $M_s$, and replace it with the solid torus $X$ constructed in Section \ref{sec:farey} (with the same dihedral angles as $T$). By Proposition \ref{prop:rigidity}, $X$ is isometric to the closure of $\Delta\cup\Delta'$, so after replacement we obtain a geometric ideal triangulation $\mathcal{D}_s$ of the filling $M_s$ (as in Theorem \ref{thm:main}). In the remainder of Section \ref{sec:voronoi}, we check that $\mathcal{D}_s$ is Delaunay.

%%%%%%%%%%%%%%%%%%%%%%%%%%

\subsection{Minkowski space} \label{sec:minkowski}
Our pictures (e.g. of the cusp link in Figure \ref{fig:hexa}) are drawn in the upper half-space model of $\mathbb{H}^3$, but we will check geometric canonicity through a computation in the Minkowski space model. This section is only a quick reminder of the formulas relating the two models, and of Epstein--Penner's convex hull construction. 

Endow ${\mathbb R}^4=\mathbb{R}^{3+1}$ with the Lorentzian product
$\left \langle (x,y,z,t) | (x',y',z',t') \right \rangle :=
xx'+yy'+zz'-tt'$. Define $$\mathcal{X}:=\{v=(x,y,z,t)\in{\mathbb R}^4 ~|~ t>0
\text{ and } \langle v | v \rangle = -1 \}.$$ Then
$\langle.|.\rangle$ restricts to a Riemannian metric on $\mathcal{X}$ and
there is an isometry $\mathcal{X} \simeq {\mathbb H}^3$, with
$\text{Isom}^+(\mathcal{X})$ a component of $SO_{3,1}({\mathbb R})$. We will
identify the point $(x,y,z,t)$ of $\mathcal{X}$ with the point at Euclidean
height $\frac{1}{t+z}$ above the complex number $\frac{x+iy}{t+z}$
in the Poincar\'e upper half-space model. Under this convention, the
closed horoball $H_{d,\zeta}$ of Euclidean diameter $d$ centered at
$\zeta=\xi+i\eta\in {\mathbb C}$ in the half-space model corresponds
to  $\{v\in \mathcal{X} ~|~ \langle v | v_{d,\zeta}\rangle\geq -1\}$, where
\begin{equation} \label{eq:isotropic}
v_{d,\zeta}=\frac{1}{d}\left (2\xi,2\eta,1-|\zeta|^2,1+|\zeta|^2\right ).
\end{equation} 
We therefore identify the horoball $H_{d,\zeta}$ with the point $v_{d,\zeta}$ of the isotropic cone (check $\langle v_{d,\zeta}|v_{d,\zeta}\rangle=0$). Similarly, the closed horoball $H_{h,\infty}$ of
points at Euclidean height no less than $h$ in the half-space model
corresponds to $\{v\in \mathcal{X} ~|~ \langle v | v_{h,\infty} \rangle \geq
-1\}$ where $v_{h,\infty}=(0,0,-h,h)$, so we identify $H_{h,\infty}$
with $v_{h,\infty}$.

Consider the following objects: a complete, oriented, cusped, finite--volume hyperbolic $3$-manifold $M$, a horoball neighborhood $H_c$ of each cusp $c$, a universal covering $\pi:{\mathbb H}^3\rightarrow M$, and the group $\Gamma
\subset \text{Isom}^+({\mathbb H}^3) \subset SO_{3,1}(\mathbb{R})$ of deck transformations of $\pi$. The $H_c$ lift to an infinite family of horoballs $(H_i)_{i\in I}$ in ${\mathbb H}^3$, corresponding to a family of isotropic vectors $(v_i)_{i\in I}$ in Minkowski space, by the above construction. The closed convex hull $C$ of $\{v_i\}_{i\in I}$ in $\mathbb{R}^{3+1}$ is $\Gamma$-invariant, and its boundary $\partial C$ comes with a natural decomposition $\widetilde{\mathcal{D}}$ into polyhedral cells. In \cite{epstein-penner, akiyoshi}, Epstein, Penner and Akiyoshi proved 
\begin{proposition}
The simplicial complex $\widetilde{\mathcal{D}}$ defines a decomposition $\mathcal{D}$ of $M$ into convex ideal hyperbolic polyhedra, by projection of each face of $\widetilde{\mathcal{D}}$ to $\mathcal{X}\simeq\mathbb{H}^3$ (with respect to $0\in\mathbb{R}^{3+1}$), and thence to $M$. The decomposition $\mathcal{D}$ of $M$ is dual to the Ford--Voronoi domain; $\mathcal{D}$ depends only on the mutual volume ratios of the $H_c$, but only a finite number of decompositions $\mathcal{D}$ arise as these volume ratios vary. \qed
\end{proposition}

Conversely, given a decomposition $\mathcal{D}$ of the manifold $M$ (still endowed with the cusp neighborhoods $H_c$) into ideal polyhedra with vertices in the cusps, in order to prove that $\mathcal{D}$ is the Epstein--Penner decomposition, we only need to consider the decomposition $\widehat{\mathcal{D}}:=\pi^*(\mathcal{D})$ of $\mathbb{H}^3$ with vertex set the centers of the horoballs $\{H_i\}_{i\in I}$, lift $\widehat{\mathcal{D}}$ to an infinite simplicial complex $\widetilde{\mathcal{D}}$ in Minkowski space $\mathbb{R}^{3+1}$ (the vertices $\{v_i\}_{i\in I}$ of $\widetilde{\mathcal{D}}$ lying over the $H_i$ in the isotropic cone, and the faces of $\widetilde{\mathcal{D}}$ being affine polyhedra), and show that $\widetilde{\mathcal{D}}$ is locally convex at each dimension--$2$ face: indeed, the projection with respect to the origin provides a homeomorphism between $\mathcal{X}\simeq \mathbb{H}^3$ and $\mathcal{D} \smallsetminus \{v_i\}_{i\in I}$; the disjoint union $\bigcup_{t\geq 1}t\widetilde{D}$ is then automatically a convex body, and its faces are exactly the cells of $\widetilde{D}$.
%is then automatically the convex hull of its own vertex set. 
In that case, we call $\mathcal{D}$ \emph{geometrically canonical}.

\begin{proposition} The codimension--one simplicial complex $\widetilde{\mathcal{D}} \subset \mathbb{R}^{3+1}$, defined by a decomposition of $M$ into polyhedra, is locally convex if and only if for every $2$--dimensional facet $F=A_1\dots A_{\sigma}$ of $\widetilde{\mathcal{D}}$ (a planar polygon in $\mathbb{R}^{3+1}$), there exists a vertex $P\notin F$ of a $3$--dimensional face of $\widetilde{\mathcal{D}}$ containing $F$, and a vertex $Q \notin F$ of the other $3$--dimensional face of $\widetilde{\mathcal{D}}$ containing $F$, such that an identity of the form 
\begin{equation}\lambda P + (1-\lambda) Q = \sum_{i=1}^{\sigma} \alpha_i A_i
~\text{ where }~ \lambda\in (0,1) ~\text{ and }~ \sum_{i=1}^{\sigma} \alpha_i >1 
\label{eq:minkonvex} \end{equation}
holds (some $\alpha_i$'s can be negative, however).
\label{obs:convex}
\end{proposition}
\begin{proof}
A more geometric way of stating the identity is as follows: if the hyperplane $\Pi \simeq \mathbb{R}^3$ is the linear span of the $A_i$'s, then the affine span of the $A_i$'s separates (in $\Pi$) the origin from the intersection of $\Pi$ with the segment $PQ$. This clearly expresses local convexity at the facet $A_1\dots A_{\sigma}$, since $P$ and $Q$ are always on opposite sides of $\Pi$ (indeed their projections to $\partial_{\infty}\mathbb{H}^3\simeq \mathbb{S}^2$ are on opposite sides of the projection of $\Pi$ to $\mathbb{H}^3$ which is a plane).  We express (\ref{eq:minkonvex}) by saying that $A_1\dots A_{\sigma}$ lies \emph{below} $PQ$ (as seen from the origin).
\end{proof}

%If $\Delta$ is an ideal tetrahedron of $\mathbb{H}^3$, and at each vertex of $\Delta$ is centered a horoball, we call $\tau_{\Delta}$ the (affine) convex hull in $\mathbb{R}^{3+1}$ of the four isotropic vectors corresponding to these four horoballs.

\subsection{Proving convexity in $\mathbb{R}^{3+1}$}
\label{sec:voronoi-final}

We now return to the ideal triangulation $\mathcal{D}_s$ of our Dehn filling, with the solid torus $X=\Delta_{1}\cup\dots\cup\Delta_{N-1}\subset \mathcal{D}_s$. For each (triangular) face $F$ of $\mathcal{D}_s$ we must prove the convexity inequality (\ref{eq:minkonvex}) of Proposition \ref{obs:convex} (applied to adjacent \emph{tetrahedra} only, hence $\sigma=3$).

If $F$ does not belong to $X$, recall that cross--ratios of tetrahedra outside $X$ in the filling $\mathcal{D}_s$ are close to what they were before filling in $\mathcal{D}$, while the volumes of the (remaining) cusp neighborhoods in the filled manifold $M_s$ are the same as in the unfilled manifold $M$: therefore, the convexity inequality (\ref{eq:minkonvex}) in $\mathcal{D}_s$, for large enough $s$, just follows from the analoguous inequality in $\mathcal{D}$.

If $F$ is one of the two faces of $\partial X$, the inequality in $\mathcal{D}_s$ again follows from the geometric canonicity of $\mathcal{D}$. Indeed, check first that the two faces of $X$ are not glued to one another: if they were (by an orientation--reversing isometry), then the sum of angles around one of the three edges of $\partial X$ would be less than or equal to $\pi$. Therefore, the face $F$ separates a tetrahedron of $X$ from a tetrahedron outside $X$. Next, consider a cover $\pi:\mathbb{H}^3\rightarrow M$ sending infinity to $c_k$ (in the upper half--space model), and the induced decomposition $\widehat{\mathcal{D}}:=\pi^*(\mathcal{D})$ of $\mathbb{H}^3$ into ideal tetrahedra. Consider a tetrahedron $\infty ABC$ of $\widehat{\mathcal{D}}$, and the neighboring tetrahedron $ABCD$ (where $A,B,C,D \in \mathbb{C}$ and $ABC$ is an acute triangle). Define $A':=B+C-A$, the symmetric of $A$ with respect to the midpoint of $B$ and $C$. Recall the tetrahedra of the solid torus $X$ are obtained by successive diagonal exchanges, beginning at the ideal triangulation of $\partial X$. Therefore, up to a permutation of $A,B,C$, the neighbor across $ABC$ of the tetrahedron corresponding (combinatorially) to $ABCD$ in $\mathcal{D}_s$, is the tetrahedron corresponding (combinatorially) to $ABCA'$. Recall the infinite simplicial complex $\widetilde{\mathcal{D}} \subset \mathbb{R}^{3+1}$. If $a,b,c,d,a',f \in \mathbb{R}^{3+1}$ are the isotropic vectors lying above the horoballs centered at $A,B,C,D,A',\infty$ (respectively), then $abcf$ and $abcd$ are neighboring faces of $\widetilde{\mathcal{D}}$ (in particular, $abc$ lies \emph{below} the segment $fd$ as seen from the origin). But by convexity of $\widetilde{\mathcal{D}}$, the facet $abc$ of $\widetilde{\mathcal{D}}$ also lies below \emph{any} segment between vertices of $\widetilde{\mathcal{D}}$, provided this segment intersects the linear span of $a,b,c$. In particular, $abc$ lies below $a'd$ (because $A',D$ lie on opposite sides of the hyperbolic plane through $A,B,C$). This is still true for the lift $\widetilde{\mathcal{D}_s}$ of the filled triangulation $\mathcal{D}_s$ if the filling slope $s$ is large enough, because the cross--ratios in $\mathcal{D}_s$ are close to those in $\mathcal{D}$. Local convexity at the face $F=ABC$ of $\mathcal{D}_s$ is proved.

\medskip

The only cases remaining are those when $F$ is an interior face of the solid torus $X$. We postpone to the end of the section the (easier) case of the ``last'' face, along which $\Delta_{N-1}$ is glued to itself, and focus on the other faces inside $X$.

Consider adjacent ideal tetrahedra $\Delta, \Delta'$ in ${\mathbb H}^3$ which are lifts from tetrahedra of the manifold $M_s$ that are consecutive tetrahedra $\Delta_{i}$ and $\Delta_{i+1}$ of the filling solid torus. We must prove that the dihedron in $\mathbb{R}^{3+1}$ between the lifts of $\Delta,\Delta'$ points ``downward'', using the criterion of Proposition \ref{obs:convex}. 

\begin{figure} [h!]\centering 
\psfrag{a}{$a$}
\psfrag{b}{$b$}
\psfrag{c}{$c$}
\psfrag{d}{$d$}
\psfrag{e}{$e$}
\psfrag{f}{$\infty$}
\psfrag{L}{$L$}
\psfrag{F}{Farey graph:}
\psfrag{q}{equator}
\psfrag{p}{$\varphi^{\nu_a}$}
\psfrag{Dm}{$\Delta_{i}$}
\psfrag{Dp}{$\Delta_{i+1}$}
\includegraphics{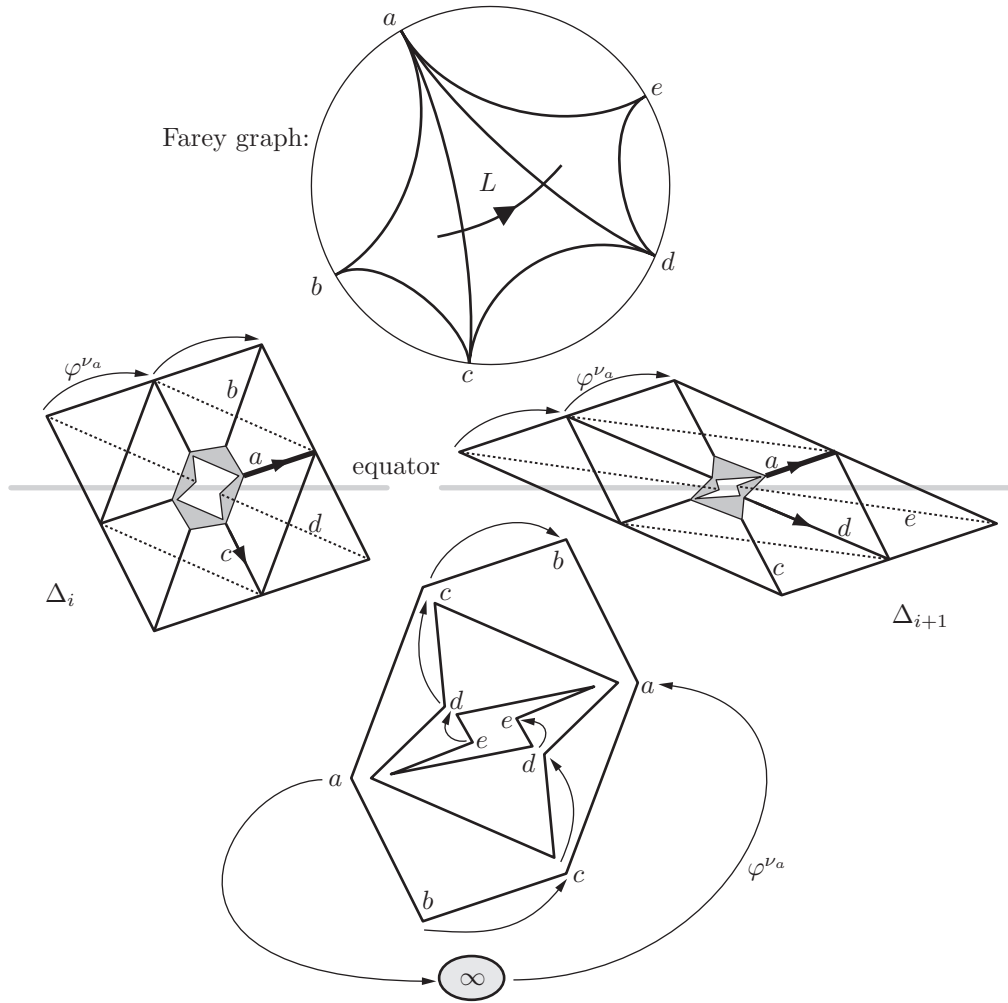}
\caption{A ``Left'' in the Farey graph corresponds to a left--handed power of $\varphi$.}
\label{fig:screwdriver}
\end{figure}

We will assume that the letter $\Omega_i$ on the Farey triangle $T_i$ is an $L$ and proceed to a careful description of the cusp link, in Figure \ref{fig:screwdriver}. Let us describe the figure.

\begin{itemize}
\item The top panel of Figure \ref{fig:screwdriver} shows a portion of the Farey graph; we name the Farey vertices $a,b,c,d,e$ so that $T_{i-1}=abc$, $T_i=acd$, $T_{i+1}=ade$ (enumerating the vertices of each triangle counterclockwise).
\item The left (resp. right) panel shows four adjacent lifts of the ideal tetrahedron $\Delta_{i}$ (resp. $\Delta_{i+1}$) in $\mathbb{H}^3$. The vertices are ideal. The direction of the equator of $\mathbb{S}^2\simeq \partial_{\infty}\mathbb{H}^3$ is materialized by a grey line. The directions $a,b,c,d,e$ of some of the ideal edges are shown. The tetrahedra in the right panel lie glued behind the tetrahedra in the left panel; the triangulation in front of the right panel thus agrees with the triangulation in the back of the left panel. In each panel, the central ideal vertex $v$, assumed to lie on the equator, has been blown up (or truncated) to depict its link, which consists of four similar Euclidean triangles drawn in grey.
\item The bottom panel puts these two ideal links together in one diagram consisting of three nested hexagons (we artificially draw each hexagon a tiny bit apart from the next one, even though they share four vertices). Each vertex of this figure corresponds to an ideal edge issued from $v$, and is marked with the slope ($a$, $b$, $c$, $d$ or $e$) of that ideal edge. (Also compare these labels with the first panel of Figure \ref{fig:hexa} page \pageref{fig:hexa}.) The four triangles between two consecutive hexagons have the same triple of angles.
\item The bottom panel also represents, up to a similarity, the endpoints in $\mathbb{C}$ of ideal edges whose other endpoint is $\infty$ in the upper half--space model of $\mathbb{H}^3$ (the point $\infty$ corresponds to the central, blown--up vertex $v$ of the previous two panels). Each triangle of the bottom panel is the vertical projection to $\mathbb{C}$ of an ideal triangle of $\mathbb{H}^3$ which, once coned off to $\infty$, yields a tetrahedron of $\mathbb{H}^3$ isometric to $\Delta_{i}$ (outer triangles) or $\Delta_{i+1}$ (inner triangles).
\item In the left (resp. right) panel we have decorated edges of slope $a$ and $c$ (resp. $a$ and $d$) with arrows. In the notation of Proposition \ref{prop:handednesses}, the loxodromy $\varphi^{\nu_a}$ is \emph{left--handed} (because $\Omega_i=L$). In these two panels, $\varphi^{\nu_a}$ acts by sending the central vertex $v$ (tail of the edge marked $a$) to the head of the edge marked $a$, and by translating all other vertices along the same direction: for example, the head of the edge marked $c$ goes to the head of the edge marked $d$.
\item This last observation allows us to understand the action of $\varphi^{\nu_a}$ on the Riemann sphere $\mathbb{C}\cup \{\infty\}$: in the bottom panel, where $v$ has been sent to $\infty$, the arrows indicate how $\varphi^{\nu_a}$ acts on the vertices of the Euclidean triangle (and $\infty$). For example, $\infty$ goes to a vertex marked $a$ and the bottom--most vertex marked $c$ goes to a vertex marked $d$. In the sequel, we must make sense of the left--handedness of this loxodromic action.
\end{itemize}

%%%%%%%%%%%%%%%%%%%%%%%%%

\begin{figure} [h!]\centering 
\psfrag{m1}{$(-1)$}
\psfrag{p1}{$(1)$}
\psfrag{z}{$(\zeta)$}
\psfrag{Z}{$(\zeta')$}
\psfrag{mz}{$(-\zeta)$}
\psfrag{mZ}{$(-\zeta')$}
\psfrag{a}{$\overrightarrow{a}$}
\psfrag{b}{$\overrightarrow{b}$}
\psfrag{c}{$\overrightarrow{c}$}
\psfrag{def}{$\begin{array}{c} 
\overrightarrow{a} = a  \, \text{exp}\,(A\sqrt{-1})\\
\overrightarrow{b} = b  \, \text{exp}\,(B\sqrt{-1})\\
\overrightarrow{c} = c  \, \text{exp}\,(C\sqrt{-1}) \end{array}$}
%\psfrag{ds}{$\overrightarrow{d_s}$}
%\psfrag{dS}{$\overrightarrow{d_{s'}}$}
%\psfrag{xi}{$x_i$}
%\psfrag{yi}{$y_i$}
%\psfrag{zi}{$\pi-w_i$}
%\psfrag{xim}{$x_{i-1}$}
%\psfrag{yim}{$y_{i-1}$}
%\psfrag{zim}{$\pi-w_{i-1}$}
\includegraphics{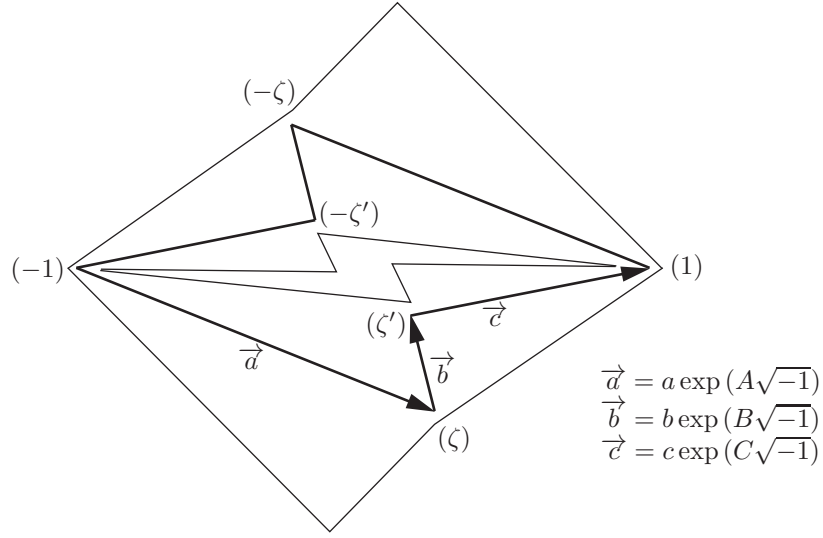}
\caption{Adjacent tetrahedra $\Delta_{i}$ and $\Delta_{i+1}$ (cusp view).}\label{fig:minkowski} \end{figure}

In order to shift to the ``Minkowski space'' aspect, we must take yet a closer look at the geometry of the link of the cusp (the following argument is taken from \cite{qf}).
In the link of the cusp, up to a complex similarity, the link of the pleated surface $\tau_i$ 
%\marginpar{$T_i$=Farey triangle, $\tau_i$=pleated surface} 
between $\Delta_{i}$ and $\Delta_{i+1}$ is
the centrally--symmetric hexagon $(-1,\zeta,\zeta',1,-\zeta,-\zeta')$ in ${\mathbb C}$, as in Figure
\ref{fig:minkowski} (which reproduces the bottom panel of Figure \ref{fig:screwdriver}): we assume that the vertices $-1,1$ both belong to the
base segments of the Euclidean triangles just inside and just outside
the hexagon. 

Let us introduce the notation
\begin{eqnarray*}\zeta+1&=&\overrightarrow{a}=a\,e^{iA} \\ \zeta'-\zeta &
=&\overrightarrow{b}=b\,e^{iB} \\1-\zeta'&=&\overrightarrow{c}=c\,e^{iC}
\end{eqnarray*} where $a,b,c\in \mathbb{R}^{>0}$ (so far $A, B, C$
are only defined modulo $2\pi$). 
The map $f:=\varphi^{\nu_a}$ now satisfies
%Above this hexagon in the half--space model lives a lift of the punctured torus $\tau_i$ with faces $(-1,\zeta,\infty),(\zeta,\infty,\zeta'),(\infty, \zeta',1)$, and which therefore admits as a deck transformation the M\"obius map $f$ that satisfies 
$f(-1)=\infty$ ; $f(\infty)=1$ ; $f(\zeta)=\zeta'$: namely,
$$f: u\mapsto 1+\frac{(\zeta+1)(\zeta'-1)}{u+1} = 1+\frac{\overrightarrow{a}\,\overrightarrow{c}}{u+1}~.$$ 
Therefore, using the notation $H_{\text{diameter, center}}$ for the horoballs of the upper half--space model (as in Section \ref{sec:minkowski}), we have $f(H_{1,\infty})=H_{|\overrightarrow{a}\overrightarrow{c}|,f(\infty)}=H_{ac,1}$. In other words, the Euclidean diameter of the horoball centered at the vertex $1$ of the hexagon is $ac$, the product of the lengths of the adjacent edges of the hexagon. 
(By an easy argument, this relationship persists if the hexagon is scaled up or down, as long as the horoball centered at infinity is $H_{1,\infty}$.)
For the same reason, the following horoballs are all sent to one another by deck transformations (in fact, by appropriate powers of $\varphi$): 
\begin{equation} \label{horoballradii}
H_{1,\infty}~;~H_{ac,-1}~;~ H_{ab,\zeta}~;~ H_{bc,\zeta'}~;~ H_{ac,1}~.
\end{equation} 
If $\zeta=\xi+\eta \sqrt{-1}$ and $\zeta'=\xi'+\eta'\sqrt{-1}$, the isotropic vectors in Minkowski space corresponding to these horoballs are respectively, using (\ref{eq:isotropic}):
\begin{equation} \begin{array}{lcccccccl} \label{lightlike}
v_{\infty}&=&&(&0,&0,&-1,&1&) \\ v_{-1}&=&\frac{1}{ac}&(&-2,&0,&0,&2&)\\
v_{\zeta}&=&\frac{1}{ab}&(&2\xi,&2\eta,&1-|\zeta|^2,&1+|\zeta|^2&)\\
v_{\zeta'}&=&\frac{1}{bc}&(&2\xi',&2\eta',&1-|\zeta'|^2,&1+|\zeta'|^2&) \\
v_{1}&=&\frac{1}{ac}&(&2,&0,&0,&2&). \end{array}\end{equation}

By Proposition \ref{obs:convex}, to prove that the dihedron at the codimension--two face (in $\mathbb{R}^{3+1}$) projecting to $(\zeta \zeta' \infty)$ is convex, it is enough to show that if $\alpha v_{\zeta}+\beta v_{\zeta'}+\gamma v_{\infty}=\lambda v_1+(1-\lambda)v_{-1}$ then $\alpha +\beta +\gamma >1$ (moreover, this will in fact take care of \emph{both} faces along which $\Delta_{i}$ touches $\Delta_{i+1}$ in the filling solid torus $X$). One easily finds the unique solution
$$\alpha=\frac{b\eta'}{c(\eta'-\eta)}~;~\beta=\frac{-b\eta}{a(\eta'-\eta)}~;~\gamma= \frac{\eta'(1-|\zeta|^2)-\eta(1-|\zeta'|^2)}{ac(\eta'-\eta)}$$
(we will not need the value of $\lambda$), hence 
$$\alpha+\beta+\gamma=1+\frac{Z}{ac(\eta'-\eta)}~~\text{ where }
Z=ab\eta'-bc\eta+\eta'(1-|\zeta|^2)-\eta(1-|\zeta'|^2)-ac(\eta'-\eta).$$ 
Observe that $\eta'>\eta$ because the triangles $-1\zeta\zeta'$ and $1\zeta'\zeta$ are counterclockwise oriented. So it is enough to prove that $Z>0$. Endow ${\mathbb C}\simeq {\mathbb R}^2$ with the usual scalar product, denoted 
``$\scalarproduct$'' to avoid confusion with scalar multiplication,
and observe that $1-|\zeta|^2=\overrightarrow{a}\scalarproduct (\overrightarrow{b}+\overrightarrow{c})$ and
$1-|\zeta'|^2=(\overrightarrow{a}+\overrightarrow{b})\scalarproduct\overrightarrow{c}$.
Hence 
\begin{eqnarray*} Z&=&
\eta'(ab+\overrightarrow{a}\scalarproduct\overrightarrow{b})
-\eta(bc+\overrightarrow{b}\scalarproduct\overrightarrow{c})
-(\eta'-\eta)(ac-\overrightarrow{a}\scalarproduct \overrightarrow{c}) \\ 
&=& abc\left [\frac{\eta'}{c}(1+\cos(A-B))-\frac{\eta}{a}(1+\cos(B-C))
-\frac{\eta'-\eta}{b}(1-\cos(A-C))\right ] \\ &=& 
-abc[\sin C(1+\cos(A\!-\!B))+\sin A(1+\cos(B\!-\!C))+\sin B(1-\cos(A\!-\!C))] \\ &=& -4abc\,\sin\frac{A+C}{2}\cos\frac{B-A}{2}\cos\frac{B-C}{2}
\end{eqnarray*} by standard trigonometric formulae. Observe that the
last expression is a well--defined function of $A,B,C \in \mathbb{R}/2\pi\mathbb{Z}$ (although each factor is defined only up to sign). Next, however, we will be careful which representatives of $A,B,C$ in $\mathbb{R}$ we pick. First, we choose for $B$ the smallest positive representative. Since the triangles $-1\zeta\zeta'$ and $1\zeta'\zeta$ are counterclockwise oriented, it follows that $B\in (0,\pi)$ and we can pick $A,C$ in $(B-\pi,B)$. Since $\overrightarrow{a}+\overrightarrow{b}+\overrightarrow{c}=2$ must also have an argument in $(B-\pi,B)$, one necessarily has
\begin{equation}
\label{eq:representatives}
-\pi<\text{min}\, \{A,C\}<0<B<\pi \hspace{8pt} \text{and} \hspace{8pt} A,C\in(B-\pi,B).
\end{equation}
In particular, to prove that $Z>0$, it only remains to show that 
\begin{equation}\label{handedness1} -\pi<\frac{A+C}{2}<0~. \end{equation}

For the deck transformation $f: u\mapsto 1+\frac{\overrightarrow{a} \overrightarrow{c}}{u+1}$ studied above, Definition \ref{def:handedness} yields ${\sf hand}\,(f)=\frac{4}{\overrightarrow{a} \overrightarrow{c}}$. But $f$ is left--handed by Proposition \ref{prop:handednesses}, so $\text{Im}\,(\overrightarrow{a} \overrightarrow{c})<0$ i.e. $A+C \in (-\pi,0)+2\pi\mathbb{Z}$. By (\ref{eq:representatives}), we have $-2\pi<A+C<\pi$ \emph{a priori}, hence in fact $-\pi<A+C<0$. Therefore (\ref{handedness1}) must hold. Geometric canonicity at the interface of tetrahedra $\Delta_{i}$ and $\Delta_{i+1}$ is proved (the argument is similar if the Farey triangle $T_i$ carries an $R$ instead of an $L$).

\smallskip

It remains to prove geometric canonicity at the core of the filling solid torus itself, where the last tetrahedron $\Delta_{N-1}$ is glued to itself along an ideal triangle. The ``hexagon'' $H_{N-1}$ of $\mathbb{C}$ has two opposite interior angles equal to $0$ and is therefore collapsed to a broken line of three segments. In (\ref{lightlike}) (and Figure \ref{fig:minkowski}), this simply translates as the identity $\zeta'=-1$; the collapsed hexagon is the broken line $(\zeta,-1,1,-\zeta)$. The radii of the horoballs centered at these vertices are computed exactly as in (\ref{horoballradii}), under the extra assumption $\zeta'=-1$.

The tetrahedra with ideal vertices $(\infty,1,-1,\zeta)$ and $(\infty,1,-1,-\zeta)$ are glued along the face $(\infty,1,-1)$, and the isotropic vectors in Minkowski space corresponding to their vertices are, following (\ref{lightlike}):

\begin{equation} \begin{array}{lcccccccl} 
v_{\infty}&=&&(&0,&0,&-1,&1&) \\ 
v_{1}&=&\frac{1}{2|1+\zeta|}&(&2,&0,&0,&2&) \\
v_{-1}&=&\frac{1}{2|1+\zeta|}&(&-2,&0,&0,&2&)\\
v_{\zeta}&=&\frac{1}{|\zeta+1|^2}&(&2\xi,&2\eta,&1-\xi^2-\eta^2,&1+\xi^2+\eta^2&)\\
v_{-\zeta}&=&\frac{1}{|\zeta+1|^2}&(&-2\xi,&-2\eta,&1-\xi^2-\eta^2,&1+\xi^2+\eta^2&)~. 
\end{array}\end{equation}

The equation $\lambda v_{\zeta}+ (1-\lambda)v_{-\zeta}=\alpha v_{\infty}+\beta v_1+\gamma v_{-1}$ has a unique solution, namely $\lambda=1/2$ and $$\alpha=\frac{|\zeta|^2-1}{|\zeta+1|^2} ~\text{ and }~ \beta=\gamma=\frac{1}{|\zeta+1|.}$$

Clearly, one will have $\alpha+\beta+\gamma>1$ if and only if $|\zeta|^2-1+2|\zeta+1|>|\zeta+1|^2$, or equivalently, $|\zeta|^2>(|\zeta+1|-1)^2$: but this relationship follows from the triangular inequality in the Euclidean triangle $(0,-1,\zeta)$. Therefore, by Proposition \ref{obs:convex}, the convexity inequality in Minkowski space is satisfied. Theorem \ref{thm:main} is proved.

\subsection{Filling on several cusps}
An analogue of Theorem \ref{thm:main} holds when several cusps undergo Dehn filling. The genericity assumptions (I--II), however, must be suitably extended.

Let $M$ be a complete hyperbolic $3$--manifold with cusps $c_1,\dots,c_k$, endowed with horoball neighborhoods $H_1,\dots, H_k$ (where $k\geq 2$). Let $l$ be an integer, $1<l\leq k$. Make the following assumptions

\begin{enumerate}[(I)]
\item The horoball neighborhoods $H_l,\dots,H_k$ are much smaller than $H_1,\dots,H_{l-1}$;
\item The decomposition $\mathcal{D}$ (before filling) consists only of ideal \emph{tetrahedra};
\item For each integer $j$ such that $l\leq j \leq k$, there exists a unique shortest path from $H_j$ to $\bigcup_{i=1}^{l-1}H_i$ in $M$;
%\item The filling slope in the cusp $c_k$ is sufficiently complicated (i.e. a finite number of ``forbidden'' slopes are avoided).
\end{enumerate}

\begin{theorem} \label{thm:main-mult}
Under the assumptions (I--III) above, for each integer $j$ such that $l\leq j \leq k$, the canonical decomposition $\mathcal{D}$ of $M$ (before filling) contains exactly two tetrahedra $\Delta_j, \Delta'_j$ with a vertex in the cusp $c_j$; moreover, $\Delta_j$ and $\Delta'_j$ are isometric and have each exactly one vertex in $c_j$ and three vertices in $\bigcup_{i=1}^{l-1}c_i$.

Moreover, for each $l\leq j \leq k$ there exists a finite set of slopes $\mathcal{X}_j$ in the cusp $c_j$ such that for any choice of slopes $s_l,\dots,s_k$ in $c_l,\dots,c_k$ satisfying $c_j \notin \mathcal{X}_j$ for each $j$, the canonical decomposition $\mathcal{D}_s$ obtained by Dehn filling along the slopes $s_l,\dots,s_k$ is combinatorially of the form
$$\mathcal{D}_s=\left ( \mathcal{D} \smallsetminus \bigcup_{j=l}^k \{\Delta_j,\Delta'_j\}\right ) \cup \bigcup_{j=l}^k \mathcal{T}_j$$
where $\mathcal{T}_j=\Delta_1^{(j)}\cup\dots\cup \Delta_{N_j-1}^{(j)}$ is a solid torus minus one boundary point, and the combinatorial gluing of the $\Delta_i^{(j)}$ (for $j$ fixed) is dictated by the continued fraction expansion of the slope $s_j$, with respect to a basis of the first homology of the cusp $c_j$ depending only on $\mathcal{D}$.
\end{theorem}

In other words, as long as the cusp neighborhoods $H_l,\dots,H_k$ are small enough and the slopes $s_l,\dots,s_k$ are long enough, Theorem \ref{thm:main} applies ``simultaneously'' to all cusps $c_l,\dots,c_k$. The proof of Theorem \ref{thm:main} transposes without major changes to Theorem \ref{thm:main-mult}, using the multicusped version of Thurston's hyperbolic Dehn surgery theorem (see e.g. Theorem 5.8.2 and the discussion immediately following it in \cite{thurstonnotes}).
% Benedetti, Petronio ? 

\section{Fillings of the Whitehead link complement} 
\label{sec:whitehead}
%%%%%%%%%%%%%%%%%%%%%%%
%%                   %%
%%   whitehead.tex   %%
%%                   %%
%%%%%%%%%%%%%%%%%%%%%%%

In this section we describe the Delaunay decompositions of all hyperbolic Dehn fillings of one cusp of the Whitehead link complement.

\subsection{Canonical decomposition before filling}

The following facts are classical; we refer to \cite{thurstonnotes} or to Weeks' program SnapPea \cite{snappea} for proofs.

\begin{figure} [h!]\centering 
\psfrag{A}{$A$}
\psfrag{B}{$B$}
\psfrag{C}{$C$}
\psfrag{D}{$D$}
\psfrag{Ap}{$A'$}
\psfrag{Bp}{$B'$}
\includegraphics{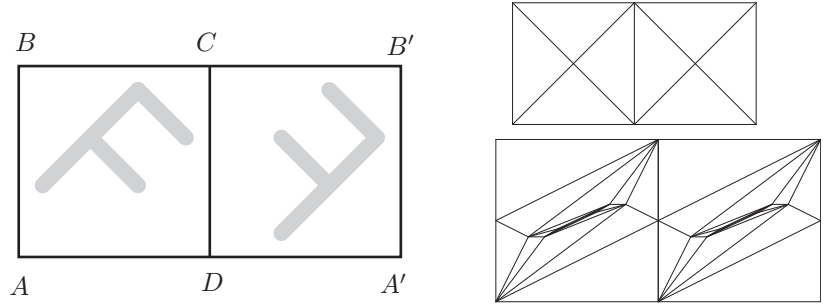}
\caption{Left: cusp view from the common tip of the square--based pyramids $Q$ and $Q'$, i.e.~from the cusp that will be filled. The $\digamma$--shaped symbol drawn on the bases of $Q$ and $Q'$ shows their identification. Right: view of the canonical decomposition from the other cusp, before (top) and after (bottom) a Dehn filling with $(k,l)=(11,8)$. In the top panel, the centers of the two squares project to the cusp that will be filled. In the bottom panel, we see that the tetrahedra in the decomposition of the filling become very close to flat, very quickly.}
\label{fig:cuspviews}
\end{figure}

Let $ABCD$ and $DCB'A'$ be two adjacent unit squares of $\mathbb{C}$ (vertices enumerated clockwise and belonging to $\mathbb{Z}[i]$, as in Figure \ref{fig:cuspviews}). Let $Q, Q'$ be the convex hulls of $\infty,A,B,C,D$ and of $\infty,D,C,B',A'$ respectively. Then $Q\cup Q'$ is a fundamental domain of the hyperbolic Whitehead link complement $M$; the face identifications are the translations of vector $\overrightarrow{AB}=i, \overrightarrow{AA'}=2$, and the hyperbolic isometry sending $A,B,C,D$ to $D,A',B',C$ respectively. Moreover, the decomposition $Q\cup Q'=M$ is the Delaunay decomposition when the horoball neighborhood of the cusp at $\infty$ is very small.

Note that $M$ has isometries that exchange the two cusps, but has no orientation--reversing isometries (so the Whitehead link is chiral).

Note also that the decomposition $Q\cup Q'$ of $M$ does not satisfy the first and second ``genericity'' assumptions of Theorem \ref{thm:main}: the cells are not tetrahedra, and the horoballs centered at $B$ and $C$, while belonging to different orbits of the stabilizer $2\mathbb{Z}\oplus i\mathbb{Z}$ of $\infty$ in the group of deck transformations, are at the same distance from the horoball at $\infty$. Thus, Theorem \ref{thm:main} does not apply directly.

\begin{proposition} \label{prop:hyperbolic-whitehead-fillings}
If $k,l$ are coprime integers, let $s$ denote the slope represented by the vector $k\overrightarrow{AA'}+l\overrightarrow{AB}$. The Dehn filling $M_s$ is hyperbolic if and only if $$\pm(k,l)\notin\{(0,1),(1,0),(1,\pm 1),(1,\pm 2)\}.$$ \end{proposition}

In the remainder of this section we assume $(k,l)$ satisfies the condition of Proposition \ref{prop:hyperbolic-whitehead-fillings} and adapt the argument of Sections \ref{sec:rivin}--\ref{sec:voronoi} to describe the Delaunay decomposition of $M_s$. This decomposition will always consist in replacing $Q\cup Q' / \langle z\mapsto z+2, z\mapsto z+i \rangle$ with a triangulated solid torus $Y$ whose exterior faces are two (triangulated) ideal quadrilaterals, which we then identify.

%%%%%%%%%%%%%%%%

\subsection{First case: $l$ is odd.} \label{sec:whitehead-odd}

If $l$ is odd, then the vector $k\overrightarrow{AA'}+l\overrightarrow{AB}=2k+il\in\mathbb{C}$ is irreducible in the lattice $\mathbb{Z}[i]$. For that reason, we can take for $Y$ the double cover of the solid torus $X$ constructed in Section \ref{sec:farey}.

More precisely, let $m \in \mathbb{P}^1\mathbb{Q}$ be the Farey vertex $\frac{l}{2k}$ (irreducible fraction). Then $m$ does not belong to $\{0,\pm 1,\pm 2,\pm \frac{1}{2},\infty\}$: the first three are ruled out because $m$ has even denominator; the last two because we assumed $\pm(k,l)\notin\{(0,1),(1,\pm 1)\}$. According to the value of $m$, choose $(p,q,r)$ as follows:
$$\begin{array}{r||c|c|c|c}
\text{if }&m<-2 &-2<m<-1 &-1<m<-1/2 &-1/2<m<0 \\
(p,q,r)=&(\infty,-1,0)&(-1,\infty,0)&(-1,0,\infty)&(0,-1,\infty)\\ \hline
\text{if }&0<m<1/2 &1/2<m<1 &1<m<2 &2<m \\
(p,q,r)=&(0,1,\infty)&(1,0,\infty)&(1,\infty,0)&(\infty,1,0)
\end{array}$$

The relative positions of $p,q,r,m$ are then exactly as in Section \ref{sec:farey}: namely, $pq$ separates $r$ from $m$; the point $m$ is not the other common Farey neighbor $r'$ of $p$ and $q$; and the line $rm$ crosses $pr'$ (not $qr'$). In particular, using the wedge notation (\ref{eq:wedge}) one has $m\wedge r \geq 3$.

Let $\theta\in (0,\pi)$ be a parameter and define 
\begin{equation} \label{eq:thetachoice}
(\theta_p,\theta_q,\theta_r):=\left \{\begin{array}{rl}
(0,\theta,\pi-\theta) & \text{if } p=\pm 1 \text{ i.e. } |m|\in (1/2,2); \\(\theta,0,\pi-\theta) & \text{if } q=\pm 1 \text{ i.e. } |m|\notin (1/2,2). 
\end{array} \right . \end{equation}
This choice will cause the ``diagonal'' edges of slope $\pm 1$ to be flat, while the edges of slope $0$ and $\infty$ will be bent.
Since $m\wedge r\geq 3$, it is straightforward to check that $(\theta_p,\theta_q,\theta_r)$ satisfies the hypothesis $(m\wedge p)\theta_p+(m\wedge q)\theta_q+(m\wedge r)\theta_r>2\pi$ of Proposition \ref{prop:pqr} if and only if $\theta$ belongs to some sub--interval $\Theta=(0,\theta_{\text{max}})\subset (0,\pi)$.

Apply now Proposition \ref{prop:pqr} and Corollary \ref{cor:completeness} with $\theta\in \Theta$. We obtain an ideal hyperbolic solid torus $X$ with dihedral angles $\theta,0,\pi-\theta$. Let $P$ be the fundamental domain of $\partial X$ defined as the ideal quadrilateral cut out by the edges of slope $0$ and $\infty$.  Let $Y$ be the double cover of $X$. Since the meridian slope is $m=\frac{l}{2k}$ and the determinant $\left | {}^1_0 {~}^{~l}_{2k} \right |$ is even, the curve of slope $\frac{1}{0}=\infty$ in $\partial X$ is homotopic to an even power of the core, and therefore lifts to a closed curve in $Y$, while the curve of slope $\frac{0}{1}=0$ does not (because $\left | {}^0_1 {~}^{~l}_{2k} \right |$ is odd). Therefore, a fundamental domain of $\partial Y$ is obtained by gluing two copies $P,P'$ of the ideal quadrilateral $P$ side by side along the edge of slope $\infty$. We view $P,P'$ as immersed in the twice--punctured torus $\partial Y$.

We now glue $P$ to $P'$ by an orientation--reversing isometry, in the same way the square bases of the pyramids $Q,Q'$ were glued together to yield the Whitehead link complement $M$ (Figure \ref{fig:cuspviews}, left). By construction, the quotient of $Y$ under this identification is homeomorphic to the Dehn filling $M_s$. The angular part of the gluing equation is automatically satisfied, since the two flat edges of $\partial Y$ (diagonals of $P,P'$) are identified, and all four non-flat edges of $\partial Y$ are identified to one edge at which the sum of dihedral angles is $\theta+(\pi-\theta)+\theta+(\pi-\theta)=2\pi$.

Therefore, the space $W$ of angle structures associated to our triangulation of $M_s$ (as in Theorem \ref{thm:rivin}) is described by setting $(\theta_p,\theta_q,\theta_r)$ as in (\ref{eq:thetachoice}) and finding all $(\theta_p,\theta_q,\theta_r)$--angle structures in the sense of Proposition \ref{prop:pqr} as $\theta$ varies freely in $\Theta \subset (0,\pi)$.

\begin{proposition} \label{prop:whitehead-odd-critical}
The volume functional has a critical point, namely a maximum, on $W$.
\end{proposition}
\begin{proof}
Exactly as in Proposition \ref{winterior}, the maximum of the (extended) volume functional is achieved at some point $z=(z_{i})_{0\leq i \leq N}$ of the closure of $W$. Using (\ref{eq:thetachoice}), the system of constraints (\ref{eq:naildown}) satisfied by $z$ now becomes
$$\begin{array}{rrrrrrr}
(&z_{0}~,& z_{1}~,& z_{2}~,&\dots~,& z_{N-1}~,&z_{N}~) \\
=~(&\pi~,& \theta~,& z_{2}~,&\dots~,& z_{N-1}~,&0~) \\
\text{or }(&\pi+\theta~,& \theta~,& z_{2}~,&\dots~,& z_{N-1}~,&0~)\end{array} $$
according to whether $|m|\in(1/2,2)$ or not.
%(where the values of $(\theta=z_{1},\dots, z_{N-1})$ remain to be chosen).

In the first case, suppose $\theta=\pi$. By the convexity condition of (\ref{racohi}), one then has $z_{0}=z_{1}=\dots=z_{h}=\pi$ where $\Delta_{h}$ is the first hinge tetrahedron. The hinge condition of (\ref{racohi}) then implies $z_{h-1}\geq z_{h}+z_{h+1}$, hence $z_{h+1}=0$. That in turn implies $z_{i}=0$ for all $i>h$ (we observed in the proof of Proposition \ref{winterior} that the sequence $(z_{i})$ is non--increasing). Therefore all tetrahedra $\Delta_{i}$ are flat, and the volume is certainly not maximal.

In the second case, suppose $\theta=\pi$. Table (\ref{tab:xiyi}) implies $\pi-\frac{z_{0}+z_{2}}{2}\geq 0$ hence $z_{2}=0$ and $z_{i}=0$ for all $i>1$: again, all $\Delta_{i}$ are flat, so the volume is certainly not maximal. 

Therefore, $\theta<\pi$. The argument of Proposition \ref{winterior} now follows through unchanged to show that no parameter $z_{i}$ for $0<i<N$ belongs to $\{0,\pi\}$. By Proposition \ref{prop:degeneracies}, all tetrahedra $\Delta_{i}$ have only positive angles (i.e.~$z\in W$).
\end{proof}

Theorem \ref{thm:rivin} applies: we have found a complete hyperbolic structure on the triangulated space $M_s$. To check that the triangulation is canonical, we only need to check the Minkowski convexity relationship (\ref{eq:minkonvex}). For interior faces of the solid torus $Y$, this is already done (Section \ref{sec:voronoi-final}). For the boundary faces, we must describe more precisely the cusp triangulation of $M_s$.

Each of the two ideal vertices of the solid torus $Y$ (projecting to the single ideal vertex of $X$) has a cusp triangulation made of nested, centrally symmetric hexagons (as in Figure \ref{fig:hexa}, right). However, by (\ref{eq:thetachoice}), two opposite angles of the outermost hexagon $H_0$ are equal to $\pi$, so the general cusp shape is a $4$--sided parallelogram. Moreover, the edges $vv',vv''$ of $H_0$ adjacent to a flat vertex $v$ have the same length: indeed, the ideal quadrilateral $\infty v'vv''$ must be a square (i.e.~its diagonals cross at a right angle), because it is a face of $Y$ and the gluing of the two isometric faces of $Y$ that yields the Dehn filling $M_s$ sends horizontal edges of one face to vertical edges of the other (e.g.~as in Figure \ref{fig:cuspviews}).

The universal cover of the cusp triangulation of $M_s$ is a union of translated copies of the cusp triangulation of $Y$. For example, up to a plane similarity, the outermost hexagons in two adjacent translates can be taken to be (for some $\zeta \in \mathbb{C} \smallsetminus \mathbb{R}$)
$$\begin{array}{ccrrrrrrc}
&(&~2\zeta-1~,& \zeta-1~,& -1~,& 1~,& \zeta+1~,& 2\zeta+1&) \\ 
\text{and } & (&-2\zeta-1~,& -\zeta-1~,& -1~,& 1~,& -\zeta+1~,& -2\zeta+1&) 
\end{array}$$
so the cusp triangles $(-1,1,\zeta+1)$ and $(-1,1,-\zeta-1)$ share an edge $(-1,1)$. We apply Proposition \ref{obs:convex} to the ideal triangle $(\infty,1,-1)$ --- by symmetry this will deal with all four triangular faces of the solid torus $Y$ (note that for proving the Minkowski convexity relationship (\ref{eq:minkonvex}), we do not care whether the two adjacent hexagons above are in the same orbit of the stabilizer of $\infty$ or not).

Following the method of Section \ref{sec:voronoi-final} (especially (\ref{horoballradii}) and the discussion that precedes it), if $\zeta=\xi+i\eta$, the isotropic vectors in $\mathbb{R}^{3+1}$ corresponding to the horoballs centered at $\infty,1,-1,\zeta+1,-\zeta-1$ are respectively
$$ \begin{array}{lcccccccl} 
v_{\infty}&=&&(&0,&0,&-1,&1&) \\ 
v_{1}&=&\frac{1}{2|\zeta|}&(&2,&0,&0,&2&) \\
v_{-1}&=&\frac{1}{2|\zeta|}&(&-2,&0,&0,&2&)\\
v_{\zeta+1}&=&\frac{1}{|\zeta|^2}&(&2\xi+2,&2\eta,&1-|\zeta+1|^2,&1+|\zeta+1|^2&)\\
v_{-\zeta-1}&=&\frac{1}{|\zeta|^2}&(&-2\xi-2,&-2\eta,&1-|\zeta+1|^2,&1+|\zeta+1|^2&)~.  \end{array} $$
The solution to $\lambda v_{\zeta+1}+ (1-\lambda) v_{-\zeta-1}=\alpha v_1+\beta v_{\infty}+\gamma v_{-1}$ satisfies $(\alpha,\beta,\gamma)=\left ( \frac{1}{|\zeta|},\frac{|\zeta+1|^2-1}{|\zeta|^2} ,\frac{1}{|\zeta|} \right )$, hence $\alpha+\beta+\gamma=1+\frac{|\zeta+1|^2-(|\zeta|-1)^2}{|\zeta|^2}>1$ according to the triangular inequality in the triangle $(0,\zeta,-1)$: by Proposition \ref{obs:convex}, the convexity inequality in Minkowski space is satisfied.

%%%%%%%%%%%%%%%%%%%%%%%%%

\subsection{Second case: $l$ is even.} \label{sec:whitehead-even}

If $l$ is even, then the vector $k\overrightarrow{AA'}+l\overrightarrow{AB}=2k+il\in\mathbb{C}$ is \emph{twice} the irreducible vector $k+i\frac{l}{2}$ in the lattice $\mathbb{Z}[i]$. For that reason, the ideal solid torus $Y$ cannot be taken to be simply a cover of $X$. Instead, we must introduce a variant of the construction of Section \ref{sec:farey}. To give a preview of the difference with Section \ref{sec:farey}, if $U\subset \mathbb{H}^3$ is a universal cover of the solid torus $Y$ we will construct below and $\langle \varphi \rangle \simeq \mathbb{Z}$ is the group of deck transformations of $U$, then for each ideal vertex $v$ of $U$, the symmetric point $v'$ of $v$ with respect to the axis of $\varphi$ is also a vertex of $U$. Moreover, $vv'$ will be an edge of the $\varphi$--invariant triangulation of $U$, and $vv'\varphi(v)\varphi(v')$ will be one of its ideal tetrahedra.

%%%%
Let $m \in \mathbb{P}^1\mathbb{Q}$ be the Farey vertex $\frac{l/2}{k}$ (reduced fraction). We have $m\notin \{\infty,0,\pm 1\}$: indeed, $\infty$ is ruled out because $m$ has odd denominator $k$ (coprime to $l$); the other possibilities are ruled out because we assumed $\pm(k,l)\notin\{(1,0),(1,\pm 2)\}$. According to the value of $m$, choose $(p,q,r)$ as in Section \ref{sec:whitehead-odd}, with the four extra possibilities
$$\begin{array}{r||c|c|c|c}
\text{if }&m=-2 &m=-1/2 &m=1/2 &m=2 \\ \hline
(p,q,r)=&(\infty,-1,0)&(0,-1,\infty)&(0,1,\infty)&(\infty,1,0)
\end{array}$$
(in fact we may switch $p,q$ in these four cases).
One then has $m\wedge r\geq 2$. Note that, unlike in Section \ref{sec:farey}, $m$ is now allowed to be the common Farey neighbor of $p$ and $q$ opposite $r$.

%Let $\theta\in (0,\pi)$ be a parameter and define $$(\theta_p,\theta_q,\theta_r):=(\theta,0,\pi-\theta).$$ This triple satisfies the hypothesis of Proposition \ref{prop:pqr} if and only if $(m\wedge r)\pi-(m\wedge q)\theta>2\pi$, which is clearly achieved for $\theta$ in some sub--interval $\Theta\subset (0,\pi)$.

Below we describe an ideal triangulation $\mathcal{D}$ for a solid torus $Y$ (with two ideal points); Proposition \ref{prop:whitehead-even-angle-structures} will then be the analogue for $\mathcal{D}$ of Proposition \ref{prop:pqr}. For convenience, we will first describe a family of tetrahedra of $\mathbb{H}^3$ whose vertices are points of $\mathbb{Z}\left [ \sqrt{-1} \right ] \subset \mathbb{P}^1\mathbb{C}\simeq \partial_{\infty}\mathbb{H}^3$, then only remember the combinatorics of the gluing of these tetrahedra. 

The sequence of Farey triangles crossed by the oriented line $\ell$ from $r$ to $m$ is $pqr=T_0,T_1,\dots,T_N=mst$ (for some Farey vertices $s,t$, and with $N\geq 1$ --- note that in Section \ref{sec:farey} we had $N\geq 2$). For every index $0\leq i \leq N$, let $x_i,y_i,z_i \in \mathbb{P}^1\mathbb{Q}$ be the vertices of $T_i$. Consider the triangulation $\mathcal{T}_i$ of $\mathbb{C}$ with vertex set $\mathbb{Z}\left [\sqrt{-1} \right ]$ and whose edges are precisely all segments of slopes $x_i,y_i,z_i$ between points of $\mathbb{Z}\left [\sqrt{-1} \right ]$. Each triangle of $\mathcal{T}_i$ is the vertical projection of an ideal triangle of $\mathbb{H}^3$ with the same triple of vertices. The union of all these ideal triangles, modulo $G:=2\mathbb{Z} \oplus \sqrt{-1}\mathbb{Z}$, is a twice--punctured torus $\tau_i$ in $\mathbb{H}^3/G$. If $0<i\leq N$ then the space between $\tau_{i-1}$ and $\tau_i$ is the union of two ideal tetrahedra $\dot{\Delta}_{i}$ and $\ddot{\Delta}_{i}$ (glued together along some of their edges). Note that the index $i=N$ is now allowed, unlike in Section \ref{sec:farey}, so that e.g.~the tetrahedron $\dot{\Delta}_{N}$ (belonging to the last pair) has an edge of slope $m$, the meridian. Also note that since $m=\frac{l/2}{k}$ and $k+\frac{l}{2}\sqrt{-1}\notin G$ (because $k$ is odd), this edge of slope $m$ runs from one of the punctures of $\tau_N$ (or $\tau_0$) to the other.

Consider now the triangulation $\{\dot{\Delta}_{i},\ddot{\Delta}_{i}\}_{1\leq i\leq N}$ as a combinatorial object only. To ``kill'' the slope $m$, we identify the edges of slope $m$ in $\dot{\Delta}_{N}$ and $\ddot{\Delta}_{N}$, and fill the remaining space with a single tetrahedron $\Delta_{N+1}$ all of whose four faces are glued to the inner faces of $\dot{\Delta}_{N} \cup \ddot{\Delta}_{N}$. This $\Delta_{N+1}$ is the tetrahedron referred to as ``$vv'\varphi(v)\varphi(v')$'' at the beginning of Section \ref{sec:whitehead-even}. We denote by $\mathcal{D}$ the triangulation $\bigcup_{i=1}^{N}\{\dot{\Delta}_{i},\ddot{\Delta}_{i}\}\cup \{\Delta_{N+1}\}$ and by $Y$ its underlying space, a twice--punctured solid torus. Note that $\mathcal{D}$ admits a combinatorial involution $\iota$ exchanging $\dot{\Delta}_{i}$ with $\ddot{\Delta}_{i}$ for all $1\leq i \leq N$ (and fixing $\Delta_{N+1}$ setwise): this $\iota$ extends the translation of $\partial Y$ that shifts one puncture to the other.

\begin{figure} [h!]\centering 
\psfrag{pb}{$\pi\!-\!b$}
\psfrag{bc}{$b\!-\!c$}
\psfrag{a}{$a$}
\psfrag{b}{$b$}
\psfrag{c}{$c$}
\psfrag{pac}{$\pi\!-\!\frac{a+c}{2}$}
\psfrag{abc}{$\frac{a-2b+c}{2}$}
\psfrag{Hn}{$H_N$}
\includegraphics{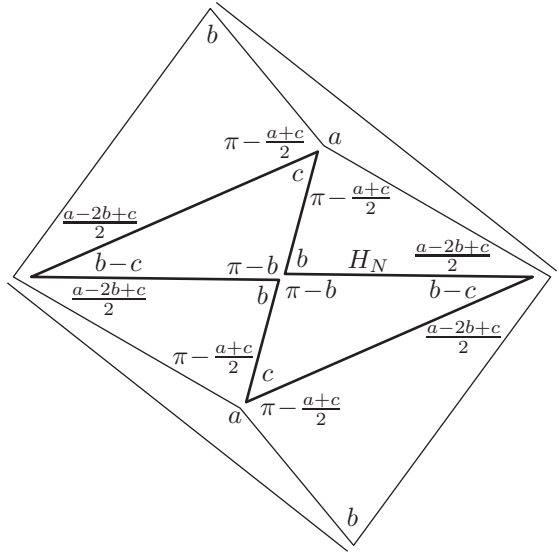}
\caption{The innermost hexagon $H_N$ along with $H_{N-1}$ and the links (Euclidean triangles) of the tetrahedra $\Delta_{N+1}, \dot{\Delta}_{N},\ddot{\Delta}_{N}$. The angles around each interior vertex sum to $2\pi$.}
\label{fig:whitehead-even}
\end{figure}

The ideal link of each of the two ideal vertices of $Y$ (which are exchanged by $\iota$) consists of nested hexagons as in Figure \ref{fig:hexa}, but the innermost hexagon is now $H_N$ (not $H_{N-1}$), and is not collapsed to a broken line of three segments. Instead, the effect of identifying the edges of slope $m$ has been to identify a pair of opposite vertices of $H_N$ (namely the inward--pointing vertices); the inside of $H_{N}$ is the union of two triangles joined by a vertex. These two triangles are two vertex links of the tetrahedron $\Delta_{N+1}$ (the other two are in the other ideal vertex of $Y$). See Figure \ref{fig:whitehead-even}.

We will not consider the full space of angle structures for our triangulation $\mathcal{D}$ of $M_s$: rather, we will restrict to $\iota$--invariant angle structures (i.e.~angle structures in which for each $1\leq i\leq N$, the tetrahedra $\dot{\Delta}_{i}$ and $\ddot{\Delta}_{i}$ have the same dihedral angles). Note that if there is an angle structure, we can always average it with its push--forward by $\iota$ to get a $\iota$--invariant angle structure.

\begin{proposition} \label{prop:whitehead-even-angle-structures}
Consider non--negative reals $\theta_p, \theta_q, \theta_r$ satisfying (\ref{eq:thetas}), namely $0<\theta_r<\pi=\theta_p+\theta_q+\theta_r$. The space of $\iota$--invariant angle structures on $\mathcal{D}$ that induce exterior dihedral angles $\theta_p, \theta_q, \theta_r$ at the edges of slope $p,q,r$ of $\partial Y$ (also called $(\theta_p,\theta_q,\theta_r)$--angle structures) is non--empty.
\end{proposition}
%\begin{remark} Proposition \ref{prop:whitehead-even-angle-structures} requires no inequality like Proposition \ref{prop:pqr}, but that is only because ``problematic'' slopes $(k,l)$ have already been ruled out. \end{remark}
\begin{proof}
As in Section \ref{sec:angle-structures}, we introduce an angle parameter $z_{i} \in (0,\pi)$ for every pair of ideal tetrahedra $\dot{\Delta}_{i},\ddot{\Delta}_{i}$ (where $1\leq i \leq N$). In what follows, $\dot{\Delta}_{i}$ and $\ddot{\Delta}_{i}$ will always be assumed to have the same dihedral angles (they are exchanged by the combinatorial symmetry $\iota$). We also denote by $z_{N+1}$ the dihedral angle of $\Delta_{N+1}$ at the edge whose slope is the only rational (Farey vertex) in $T_N\cap T_{N-1} \smallsetminus T_{N-2}$. Using these conventions and writing $(a,b,c):=(z_{N-1},z_{N},z_{N+1})$, it is easy to see that the triples of dihedral angles of the ideal tetrahedra are as follows: 
\begin{equation} \begin{array}{ccrrr} \label{eq:last-angles}
\dot{\Delta}_{N},\ddot{\Delta}_{N}&:&
(~b~,&\pi-\frac{a+c}{2}~,& \frac{a-2b+c}{2}~) \\ 
\Delta_{N+1}&:&(~c~,&\pi-b~,&b-c~) \end{array}\end{equation}
(see also Figure \ref{fig:whitehead-even}). For $1\leq i <N$, the dihedral angles of $\dot{\Delta}_{i},\ddot{\Delta}_{i}$ are simply given by Table (\ref{tab:xiyi}). In keeping with Table (\ref{tab:xiyi}), we consider $z_{N}$ to be a non--hinge parameter.

Recall that $N\geq 1$: analogously to (\ref{eq:naildown}) -- (\ref{racohi}), we are thus looking for sequences of the form 
$$\begin{array}{rrrrrrr}
(&z_{0}~,& z_{1}~,& z_{2}~,&\dots~,& z_{N}~,&z_{N+1}~) \\
=~(&\pi+\theta_q~,& \pi-\theta_r~,& z_{2}~,&\dots~,& z_{N}~,&z_{N+1}~) \end{array} $$
subject to the conditions
$$ \left \{ \begin{array}{ll}
 z_{i-1}>z_{i}+z_{i+1} 
& \text{if $z_{i}$ is a hinge parameter (\emph{hinge condition});} \\
 z_{i-1}+z_{i+1}>2z_{i} 
& \text{if 
%$z_{i}$ is a non--hinge parameter (e.g.~$i=N$)
not (\emph{convexity condition}), e.g. $i=1$ or $N$
;} \\
 0<z_{i}<\pi 
& \text{for all $2\leq i \leq N$ (\emph{range condition});} \\
 0<z_{2}<\pi-\theta_q 
& \text{as in (\ref{racohi}) above;} \\
 0<z_{N+1}<z_{N} 
& \text{which follows from (\ref{eq:last-angles}).} \\ \end{array} \right . $$
To find such a sequence, the argument that finishes Section \ref{sec:angle-structures} follows through essentially unchanged: we construct a convex positive decreasing sequence $(z_{i})_{0\leq i \leq h}$ where $h$ is the smallest hinge index (or $h=N+1$ if there are no hinges), then set e.g.~$z_{i+1}=\varepsilon z_{i}$ (inductively) for all $i\geq h$ and a fixed small $\varepsilon>0$.
\end{proof}

Finally, we must glue the faces of the solid torus $Y$ together to form the Dehn filling $M_s$ of the Whitehead link complement. This is performed exactly as in Section \ref{sec:whitehead-odd}: we set $(\theta_p,\theta_q,\theta_r)$ as in (\ref{eq:thetachoice}) for $0<\theta<\pi$, so that the faces of $\partial Y$ become two ideal quadrilaterals $P,P'$ with edges of slopes $0$ and $\infty$; then glue $P$ to $P'$ by an orientation--reversing homeomorphism sending the edges of slope $0$ of $P$ to the edges of slope $\infty$ of $P'$ (and conversely). The angular gluing equations are automatically satisfied. 

Therefore, the full space $W$ of $\iota$--invariant angle structures for $\mathcal{D}$ is obtained by letting $\theta$ range over $(0,\pi)$ and finding all $(\theta_p,\theta_q,\theta_r)$--angle structures in the sense of Proposition \ref{prop:whitehead-even-angle-structures}.

\begin{proposition}
The volume functional has a critical point, namely a maximum, on $W$.
\end{proposition}
\begin{proof}
As in Proposition \ref{winterior}, the maximum of the (extended) volume functional is achieved at some point $z$ of the closure $\overline{W}$ of $W$. Using (\ref{eq:thetachoice}), the system of constraints (\ref{eq:naildown}) becomes
$$\begin{array}{rrrrrrr}
(&z_{0}~,& z_{1}~,& z_{2}~,&\dots~,& z_{N}~,&z_{N+1}~) \\
=~(&\pi~,& \theta~,& z_{2}~,&\dots~,& z_{N}~,&z_{N+1}~) \\
\text{or }(&\pi+\theta~,& \theta~,& z_{2}~,&\dots~,& z_{N}~,&z_{N+1}~) \end{array} $$
according to the value of $m$.
%(where the values of $(\theta=z_{1},\dots, z_{N+1})$ remain to be chosen).

The assumption $\theta=\pi$ leads to a contradiction exactly as in the proof of Proposition \ref{prop:whitehead-odd-critical}. Therefore $\theta<\pi$.

By (\ref{eq:last-angles}), $\dot{\Delta}_{N}$ and $\ddot{\Delta}_{N}$ have a dihedral angle equal to $b:=z_{N}$, while $\Delta_{N+1}$ has an angle $\pi-b$. On the other hand, a tetrahedron of $\mathcal{D}$ is flat at $z\in \overline{W}$ if and only if one (and therefore all) of its angles belong to $\{0,\pi\}$ (Proposition \ref{prop:degeneracies}). Thus, $\dot{\Delta}_{N},\ddot{\Delta}_{N}$ are flat if and only if $\Delta_{N+1}$ is flat (i.e.~$b\in\{0,\pi\}$). The argument of Proposition \ref{winterior} then follows through: at $z$, if some tetrahedra were flat, all would be flat and the volume would be $0$; absurd. Thus $z\in W$.
\end{proof}

To apply Theorem \ref{thm:rivin}, we only need to make sure that the critical point (maximum) of $\mathcal{V}$ on the space $W$ of $\iota$--invariant angle structures is also critical (maximal) in the space of \emph{all} angle structures: but that is clear since by concavity of the volume functional (Fact \ref{fact:infinite-derivative}), the volume can only go up when we average an angle structure with its push--forward by $\iota$. Theorem \ref{thm:rivin} does apply: we have found a complete hyperbolic structure on the triangulated space $M_s$. To check that the triangulation is canonical, we only need to check the Minkowski convexity relationship (\ref{eq:minkonvex}).
For boundary faces of $Y$, the situation is exactly the same as in Case $1$ (odd $l$). For interior faces of $Y$ not bounding the ``extra'' tetrahedron $\Delta_{N+1}$, we proceed as in Section \ref{sec:voronoi-final}: the only new argument needed is an analogue of Proposition \ref{prop:handednesses} (predicting the handednesses of powers of the core curve of $Y$), namely

\begin{proposition} \label{prop:handednesses-whitehead-odd}
Let $T_i$ be a Farey triangle such that $0<i<N$ and let $x\in \mathbb{P}^1\mathbb{Q}$ be the Farey vertex $T_{i-1}\cap T_i\cap T_{i+1}$. Consider a properly embedded line $\gamma_x$ of slope $x$ in $\partial Y$ (running between two cusps), and a lift $\widehat{\gamma_x}$ of $\gamma_x$ to a universal cover $U\subset \mathbb{H}^3$ of $Y$ (running between two ideal points). The deck transformation of $U$ that sends the initial point of $\widehat{\gamma_x}$ to the final point is left--handed (resp. right--handed) if and only if the Farey triangle $T_i$ carries a letter $L$ (resp. $R$).
\end{proposition}
\begin{proof}
The proof is exactly as in Section \ref{sec:handedness}. The key argument that the integral $\lambda_x$ of the longitude $1$--form along $\widehat{\gamma_x}$ stays less than $\pi$ is only easier, because the ``longest'' curve $\gamma_m$ runs only around one half, not all, of the meridian of $U$ (connecting some ideal vertex to its symmetric with respect to the axis of $U$); thus $\lambda_m=\pi$ and $\lambda_x<\pi$.
\end{proof}

The only remaining case of the Minkowski convexity relationship (\ref{eq:minkonvex}) is at the faces of $\Delta_{N+1}$. According to our picture of the cusp triangulation (Figure \ref{fig:whitehead-even}), we can assume that the innermost hexagon $H_N$ has vertices at $$-1~,~0~,~\zeta~,~1~,~0~,~-\zeta$$ and look at the interface $\zeta \infty 0$ between ideal tetrahedra $1\zeta \infty 0$ and $-1 \zeta \infty 0$. 

Following the method of Section \ref{sec:voronoi-final}, if $\zeta=\xi+i\eta$, the isotropic vectors in $\mathbb{R}^{3+1}$ corresponding to the horoballs centered at $\infty,0,\zeta,1,-1$ are respectively
$$ \begin{array}{lcccccccl} 
v_{\infty}&=&&(&0,&0,&-1,&1&) \\ 
v_{0}&=&\frac{1}{|\zeta|}&(&0,&0,&1,&1&) \\
v_{\zeta}&=&\frac{1}{|\zeta||\zeta-1|}&(&2\xi,&2\eta,&1-|\zeta|^2,&1+|\zeta|^2&)\\
v_{1}&=&\frac{1}{|\zeta-1|}&(&2,&0,&0,&2&)\\
v_{-1}&=&\frac{1}{|\zeta-1|}&(&-2,&0,&0,&2&)~.  \end{array} $$
The solution to $\lambda v_{1}+ (1-\lambda) v_{-1}=\alpha v_{\infty}+\beta v_0+\gamma v_{\zeta}$ satisfies $(\alpha,\beta,\gamma)=\left ( \frac{1}{|\zeta-1|},\frac{|\zeta|}{|\zeta-1|} ,0 \right )$, hence $\alpha+\beta+\gamma=\frac{|\zeta|+1}{|\zeta-1|}>1$ according to the triangular inequality in the triangle $(0,1,\zeta)$: by Proposition \ref{obs:convex}, the convexity inequality in Minkowski space is satisfied.

\subsection{Delaunay decompositions and elementary Kleinian groups} \label{sec:elementary-extension}

\begin{remark}
If $U\subset \mathbb{H}^3$ is a (triangulated) universal cover of the solid torus $Y$ and $\langle \varphi \rangle$ is the group of deck transformations of $U$, we mentioned at the beginning of Section \ref{sec:whitehead-even} that for each ideal vertex $v$ of $U$, the symmetric $v'$ of $v$ with respect to the axis of $\varphi$ is also a vertex of $U$, and $\Delta:=vv'\varphi(v)\varphi(v')$ is an ideal tetrahedron of $U$ (projecting to $\Delta_{N+1}$). By duality between the Ford--Voronoi domain and the canonical triangulation, the last computation of Section \ref{sec:whitehead-even} amounts to checking the following (easy) fact: if all vertices of $U$ are endowed with horoballs of the same size, then the center of $\Delta$ is nearer to the horoballs centered at the vertices of $\Delta$ than to any other horoballs.
\end{remark}

More generally, if $n\geq 3$, let $G:=\langle \varphi, \psi \rangle \subset \text{Isom}^+(\mathbb{H}^3)$ be an elementary group generated by a loxodromy $\varphi$ and an order--$n$ rotation $\psi$ with the same axis $\delta$ (note that Section \ref{sec:whitehead-even} amounted to the case $n=2$, and Section \ref{sec:farey} to the case $n=1$). Let $\mathcal{O}:=Gp \subset \partial_{\infty}\mathbb{H}^3$ be a generic ideal orbit of $G$; if $h_p$ is a horoball centered at $p$, all horoballs in the $G$--orbit of $h_p$ come equally close to the line $\delta$. The convex hull of $\mathcal{O}$ projects modulo $\varphi$ to an $n$--times punctured solid torus $X$ whose boundary is pleated along a certain ideal triangulation in which all vertices have the same degree (generically $6$, exceptionally $4$; for simplicity let us assume the generic situation). The convex hull construction in Minkowski space $\mathbb{R}^{3+1}$ yields a decomposition of $X$ into ideal polyhedra with respect to the horoballs $G h_p$. The central polyhedron is the convex hull $Q$ of $\langle \psi \rangle p \cup \varphi(\langle \psi \rangle p)$, namely an ideal hyperbolic antiprism with regular $n$--sided bases (glued together \emph{via} $\varphi$): indeed, it is easy to check that the center of $Q$ is closer to the horoballs centered at the vertices of $Q$ than to any other horoballs of the $G$--orbit. 

It is possible that $Q$ is the only cell of $X$. Otherwise, we claim that the remaining cells between $Q$ and $\partial X$ are tetrahedra glued together according to diagonal exchanges and Farey--type combinatorics: namely, $\partial X/\psi$ is a once--punctured torus with ideal edges of slope $p,q,r \in \mathbb{P}^1\mathbb{Q}$ for some arbitrary marking (these slopes are mutual Farey neighbors). The meridian of X defines the $n$-th power of an irreducible element of $H_1(\partial X/\psi,\mathbb{Z})$, and therefore also a slope $m \in \mathbb{P}^1\mathbb{Q}$. Since $m$ is the slope of the base edges of the antiprism $Q$, if $Q$ is the only cell in $X$ then $m\in \{p,q,r\}$. Otherwise, we may as in Section \ref{sec:farey} assume that the Farey edge $pq$ separates $m$ from $r$ and follow the line $\ell$ from $r$ to $m$ to construct a (combinatorial) ideal decomposition $\mathcal{D}$ of $X$.

In fact, the following ``Gauss--Bonnet type'' result (left as an exercise) is a simple generalization of the method worked out in this paper. It uses the fact that the antiprism $Q$ (like any convex ideal hyperbolic polyhedron, see \cite{rivin-rigidity,fg-rigidity}) is uniquely determined up to isometry by its dihedral angles.

\begin{theorem}
Consider non--negative reals $\theta_p, \theta_q, \theta_r$ satisfying (\ref{eq:thetas}), namely $0<\theta_r<\pi=\theta_p+\theta_q+\theta_r$. There exists a hyperbolic $n$-times punctured solid torus $X$, decomposed into convex ideal polyhedra according to the combinatorics of $\mathcal{D}$ and with exterior dihedral angles $\theta_p, \theta_q,\theta_r$ at the edges of slope $p,q,r$, if and only if $$(m\wedge p)\theta_p+(m\wedge q)\theta_q+(m\wedge r)\theta_r>\frac{2\pi}{n}.$$
Moreover, $X$ is then unique up to isometry and $\mathcal{D}$ is the Delaunay decomposition of $X$. \qed
\end{theorem}

%\begin{figure}[h!] \centering
%\psfrag{0}{$0$}
%\includegraphics{xyz.eps}
%\caption{blabla \label{xxx}} \end{figure}

\begin{flushright}
% ADDRESSES
\end{flushright}

\end{document}